# Lefschetz fibrations over the disc

Nikos Apostolakis, Riccardo Piergallini and Daniele Zuddas

## Abstract

We provide a complete set of moves relating any two Lefschetz fibrations over the disc having as their total space the same four-dimensional 2-handlebody up to 2-equivalence. As a consequence, we also obtain moves relating diffeomorphic three-dimensional open books, providing a different approach to an analogous previous result by Harer.

## 1. Introduction

As Harer showed in [11], any four-dimensional 2-handlebody $W$ can be represented by a topological (achiral) Lefschetz fibration over the disc, that is, a smooth map $W \to B^2$ whose generic fibre is an orientable bounded surface and whose singularities are topologically equivalent to complex non-degenerate ones. Harer's argument is based on Kirby calculus [14].

An alternative approach to the same result was provided in [15, Remark 2.3]. This is based on the characterization of allowable Lefschetz fibrations (see Section 6) as those smooth maps that admit a factorization $W \xrightarrow{p} B^2 \times B^2 \xrightarrow{\pi} B^2$, where $p : W \to B^2 \times B^2$ is a covering simply branched over a braided surface and $\pi$ is the canonical projection. The two other ingredients of the proof are Montesinos's representation of four-dimensional 2-handlebodies as coverings of $B^4$ simply branched over ribbon surfaces [16] and Rudolph's procedure for isotoping any orientable ribbon surface to a braided surface [20].

In this paper, we use the second approach together with the branched covering interpretation of Kirby calculus given by Bobtcheva–Piergallini in [3], to relate different Lefschetz fibrations representing the same four-dimensional 2-handlebody up to 2-equivalence by means of certain moves $S, T$ and $U$ on their monodromy representation.

These monodromy moves are described in Section 7. Move $S$ (Figure 33) is nothing but the well known positive or negative Hopf stabilization, and it corresponds to adding a pair of cancelling 1- and 2-handles to the handlebody, while move $T$ (Figure 36) is new, and roughly speaking it corresponds to a 2-handle sliding. Both such moves are applied only to allowable Lefschetz fibrations (see Section 6). On the contrary, move $U$ (Figure 38) is just used to transform any Lefschetz fibration into an allowable one.

Namely, our main result is the following theorem in Section 8.

THEOREM A. *Any two allowable Lefschetz fibrations $f : W \to B^2$ and $f' : W' \to B^2$ represent 2-equivalent four-dimensional 2-handlebodies $H_f$ and $H_{f'}$ if and only if they are related by fibred equivalence and the moves $S$ and $T$. Moreover, the allowability hypothesis can be relaxed by using in addition move $U$.*

2010 *Mathematics Subject Classification* 55R55, 57N13 (primary), 57M12, 57R65 (secondary).

Partially supported by a CUNY Community College Collaborative Incentive Grant and a PSCCUNY Cycle 39 Research Award (Nikos Apostolakis). Supported by Regione Autonoma della Sardegna with funds from PO Sardegna FSE 2007–2013 and L.R. 7/2007 'Promotion of scientific research and technological innovation in Sardinia'. Also thanks to ESF for short visit grants within the program 'Contact and Symplectic Topology' (Daniele Zuddas).

Here is a very sketchy outline of the proof of Theorem A. According to [**15**], the two allowable Lefschetz fibrations are realized as simple coverings of $B^2 \times B^2$ branched over braided surfaces (see Section 6). Such braided surfaces, endowed with the labelling that encodes the monodromy of the coverings, are special cases of labelled ribbon surfaces representing 2-equivalent, four-dimensional 2-handlebodies as branched coverings of $B^4$ (see Section 5). Therefore, they can be related by a finite sequence of isotopy and covering moves (Figures 3 and 24) given in [**3**]. Then, we perform on these moves a streamlined version of the Rudolph's braiding procedure [**20**], which retracts labelled ribbon surfaces onto labelled braided surfaces (see Section 4). The result is a quite large set of moves on labelled braided surfaces, and the last part of the proof, carried out in Section 8, consists in reducing it, up to braided isotopy, to only two moves corresponding to the monodromy moves $S$ and $T$.

The same argument also gives the following theorem in Section 9. Here, the extra move $P$ (Figure 66) corresponds to making connected sum with $\mathbb{C}P^2$, whereas move $Q$ consists in adding a pair contiguous opposite Dehn twists to the monodromy sequence of the Lefschetz fibration.

THEOREM B. *Two allowable Lefschetz fibrations over $B^2$ represent four-dimensional 2-handlebodies with diffeomorphic oriented boundaries if and only if they are related by fibred equivalence, the moves $S$ and $T$ of Section* 7, *and the moves $P$ and $Q$.*

Theorems A and B can be considered as four-dimensional analogues of the equivalence theorem for three-dimensional open books proved by Harer in [**12**]. In fact, such open books naturally arise as boundary restrictions of Lefschetz fibrations. Then, by considering the boundary restrictions $\partial S, \partial T$ and $\partial P$ of the moves $S, T$ and $P$, we also derive the next theorem in Section 9. We remark that, in contrast to Harer's moves, our moves can be completely described in terms of the open book monodromy.

THEOREM C. *Two open books are supported by diffeomorphic oriented three-manifolds if and only if they are related by fibred equivalence and the moves $\partial S, \partial T$ and $\partial P$.*

The paper is organized as follows. Section 2 is devoted to ribbon surfaces and to 1-isotopy between them. Sections 3 and 4 deal with braided surfaces and the Rudolph's braiding procedure. In Section 5, we review the branched covering representation of four-dimensional 2-handlebodies and adapt the covering moves to the present aim. In Sections 6 and 7, we recall the branched covering representation of Lefschetz fibrations and define the equivalence moves for them. Finally, in Sections 8 and 9 we establish the three equivalence theorems stated above.

## 2. *Ribbon surfaces*

A regularly embedded smooth compact surface $S \subset B^4$ is called a *ribbon surface* if the Euclidean norm restricts to a Morse function on $S$ with no local maxima in Int $S$. Assuming that $S \subset R^4_+ \subset R^4_+ \cup \{\infty\} \cong B^4$, where $\cong$ stands for the standard orientation preserving conformal equivalence, this property is smoothly equivalent to the fact that the fourth Cartesian coordinate restricts to a Morse height function on $S$ with no local minima in Int $S$. Such a surface $S \subset R^4_+$ can be horizontally (preserving the height function given by the fourth coordinate) isotoped to make its orthogonal projection in $R^3$ a self-transversal immersed surface, whose double points form disjoint arcs as in Figure 1(a). We call the orthogonal projection $\pi(S) \subset R^3$ a *three-dimensional diagram* of $S$.

Actually, any immersed compact surface $S \subset R^3$ with all self-intersections as above and no closed components is the three-dimensional diagram of a ribbon surface. This can be obtained by pushing $\mathrm{Int}\, S$ inside $\mathrm{Int}\, R^4_+$ in such a way that all self-intersections disappear. Moreover, it is uniquely determined up to vertical isotopy.

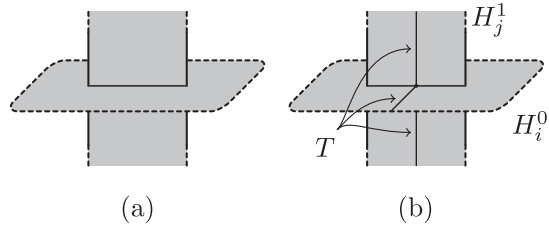

Figure 1. *Ribbon intersection.*

In the following, we will omit the projection $\pi$ and use the same notation for a ribbon surface in $B^4$ and its three-dimensional diagram in $R^3$, the distinction between them being clear from the context.

Any ribbon surface $S$ admits a handlebody decomposition with only 0- and 1-handles induced by the height function. Such a 1-handlebody decomposition $S = (H^0_1 \sqcup \cdots \sqcup H^0_m) \cup (H^1_1 \sqcup \cdots \sqcup H^1_n)$ is called *adapted*, if each ribbon self-intersection of its three-dimensional diagram involves an arc contained in the interior of a 0-handle and a proper transversal arc in a 1-handle (cf. [**21**]). Then, looking at the three-dimensional diagram, we have that the 0-handles $H^0_i$ are disjoint non-singular discs in $R^3$, while the 1-handles $H^1_j$ are non-singular bands in $R^3$ attached to the 0-handles and possibly passing through them to form ribbon intersections like the one shown in Figure 1(b). Moreover, we can think of $S$ as a smoothing of the frontier of $((H^0_1 \sqcup \cdots \sqcup H^0_m) \times [0,1]) \cup ((H^1_1 \sqcup \cdots \sqcup H^1_n) \times [0,1/2])$ in $R^4_+$.

A ribbon surface $S \subset R^4_+$ endowed with an adapted handlebody decomposition as above will be referred to as an *embedded two-dimensional 1-handlebody*.

A convenient way of representing a ribbon surface $S$ arises from the observation that its three-dimensional diagram, considered as a two-dimensional complex in $R^3$, collapses to a graph $T$. We can choose $T = \pi(P)$ for a smooth simple spine $P$ of $S$ (simple means that all the vertices have valency one or three), which intersects each 1-handle $H^1_j$ along its core. Moreover, we can also assume $T$ to meet each ribbon intersection arc of $S$ at exactly one 4-valent vertex, as in Figure 1(b) where the fourth edge of $T$ in the back is not visible. The inverse image of such a 4-valent vertex of $T$ consists of two points, in the interior of two distinct edges of $P$, while the projection restricted over the complement of all 4-valent vertices of $T$ is injective.

Therefore, $T$ has vertices of valency 1, 3 or 4. We call *singular vertices* the 4-valent vertices located at the ribbon intersections, and *flat vertices* all the other vertices. Moreover, we assume $T$ to have three distinct tangent lines at each flat 3-valent vertex and two distinct tangent lines at each singular vertex.

Up to a further horizontal isotopy of $S$, we can contract its three-dimensional diagram to a narrow regular neighbourhood of the graph $T$. Then, by considering a planar diagram of $T$, we easily get a new diagram of $S$, consisting of a number of copies of the local spots shown in Figure 2, and some non-overlapping flat bands connecting those spots. We call this a *planar diagram* of $S$.

We emphasize that a planar diagram of $S$ arises as a diagram of the *pair* $(S,T)$ and this is the right way to think about it. However, we omit to draw the graph $T$ in the pictures of planar diagrams, since it can be trivially recovered, up to diagram isotopy, as the core of the diagram itself. In particular, the diagram crossings and the singular vertices of $T$ are located at the centres of the copies of the two rightmost spots in Figure 2, while the flat vertices of it are located at the centres of copies of the two leftmost spots.

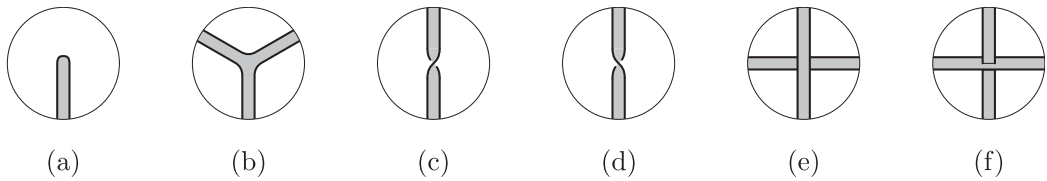

(a) (b) (c) (d) (e) (f)

Figure 2. *Local models for planar diagrams.*

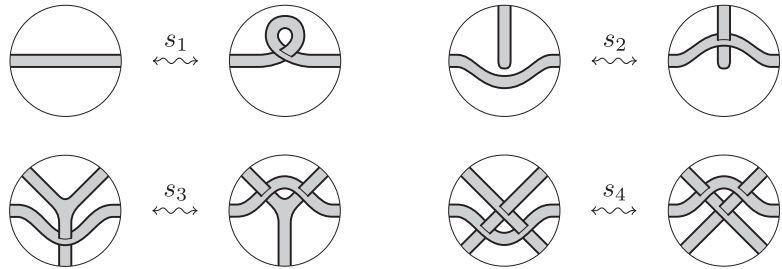

Figure 3. *1-isotopy moves for planar diagrams.*

Of course, a planar diagram determines a ribbon surface $S$ only up to vertical isotopy. Namely, the three-dimensional height function (and the four-dimensional one as well) cannot be determined from the planar diagram, apart from the obvious constrains imposed by the consistency with the local configurations of Figure 2.

Ribbon surfaces will be always represented by planar diagrams and considered up to vertical isotopy (in the sense just described above). Moreover, planar diagrams will be always considered up to planar diagram isotopy, that is ambient isotopy of the plane containing them.

Following [**3**], two ribbon surfaces $S, S' \subset R_+^4$ are said to be *1-isotopic* if there exists a smooth ambient isotopy $(h_t)_{t \in [0,1]}$ such that: (1) $h_1(S) = S'$; (2) $S_t = h_t(S)$ is a ribbon surface for every $t \in [0, 1]$; (3) the projection of $S_t$ in $R^3$ is an honest three-dimensional diagram except for a finite number of critical $t$'s. Such equivalence relation between ribbon surfaces can be interpreted as embedded 1-deformation of embedded two-dimensional 1-handlebodies, and this is the reason for calling it 1-isotopy. Whether or not 1-isotopy coincides with isotopy is unknown, but this problem is not relevant for our purposes.

All we need to know here is that two ribbon surfaces are 1-isotopic if and only if their three-dimensional diagrams are related by three-dimensional isotopies and the moves depicted in Figure 3. This has been proved in [**3**, Proposition 1.3].

We observe that the moves $s_1$ to $s_4$ are described in terms of planar diagrams. An analogous expression of three-dimensional isotopies in terms of certain moves of planar diagrams has been provided by [**4**, Proposition 10.1]. Since this aspect will be crucial in the following, we give a complete account of that result in the proof of Proposition 1.

In order to express three-dimensional isotopy of three-dimensional diagrams of ribbon surfaces in terms of planar diagrams, it is convenient to consider the special case when all ribbon intersections are terminal, that is, they appear only at the ends of the bands (and never in the middle of them, as in the rightmost spot in Figure 2).

A planar diagram with this property will be called a *special planar diagram*. Figure 4 depicts the two different local configurations that replace (f) of Figure 2, when dealing with a special planar diagram. Notice that, in the previous context, (g) and (h) can be seen as combinations of (f) and (a) of Figure 2.

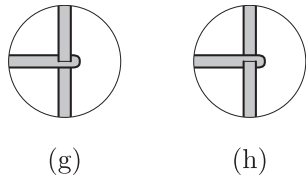

Figure 4. *Local models for ribbon intersections in special planar diagrams.*

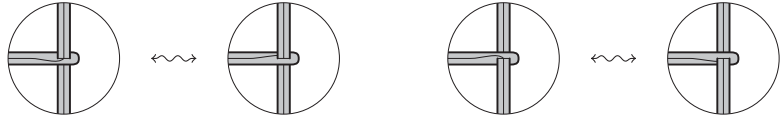

Figure 5. *Graph moves at a singular vertex.*

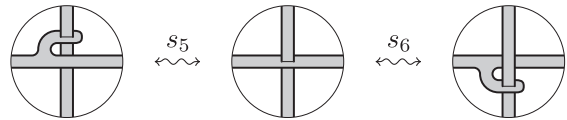

Figure 6. *Making ribbon intersections terminal.*

In this case, each singular vertex of the core graph $T \subset S$ has valency 3 and its inverse image in $P$ consists of a 1-valent vertex and a point in the interior of an edge. We still assume that $T$ has two distinct tangent lines at any singular vertex.

The three edges of $T$ converging at a ribbon intersection arc of $S$ are drawn in Figure 1(b). When we think of $T$ as a graph embedded in $S$ and represent $S$ by a planar diagram, they can be recovered from the planar diagram only up to the moves in Figure 5.

Actually, any planar diagram can be transformed into a special one, by performing one of the two moves of Figure 6 at each ribbon intersection. To get terminal ribbon intersections of type (h) instead of (g), we could also introduce symmetric moves $\bar{s}_5$ and $\bar{s}_6$, but these would derive from $s_5$ and $s_6$ in the presence of the moves described below.

Figures 7 and 8 present the three-dimensional isotopy moves for planar diagrams. They are grouped into the two figures depending on whether half-twists are involved or not.

Note that moves $s_7$ and $s_8$ do not change the topology of the planar diagram of the surface, but they change the structure of its core graph, and this is the reason why they are there. Moreover, some of the moves could be derived from the others and they are included for the sake of convenience. For example, move $s_9$ can be obtained as a combination of $s_{11}$ and $s_{13}$ modulo $s_7$, while move $s_{20}$ is a consequence of $s_7$, $s_{23}$ and $s_{24}$.

On the other hand, for each move $s_i$ in the Figures 6–8 one can consider the symmetric move $\bar{s}_i$ obtained from $s_i$ by reflection with respect to the projection plane, which reverses half-twists, crossings and ribbon intersections, interchanging local models (g) and (h). Such symmetric moves coincide with the original ones for $i = 7, \ldots, 10, 19$, while the symmetry interchanges $s_i$ and $s_{i+1}$ for $i = 11, 13, 21$. In all the other cases, the symmetry produces new moves, which nevertheless can be derived from those in Figures 6–8. In particular, moves $\bar{s}_5$ and $\bar{s}_6$, as well as the symmetric moves of Figure 7 not considered above, are generated by the original ones modulo the other moves in the same Figure 7. We leave the easy verification of this fact to the reader.

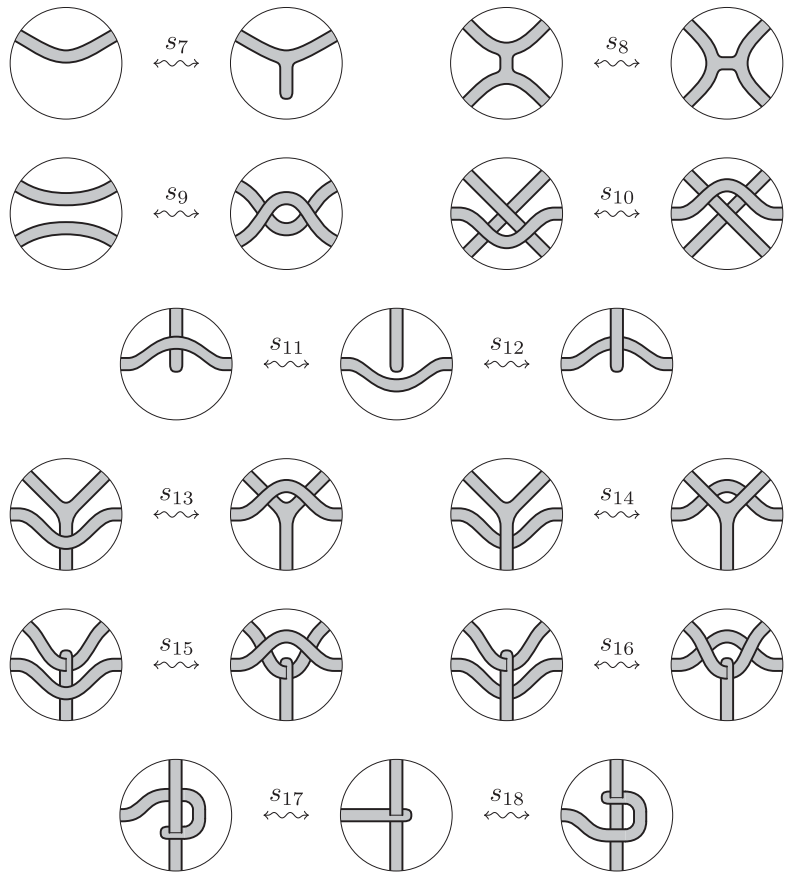

Figure 7. *Flat isotopy moves.*

PROPOSITION 1. *Two planar diagrams represent 1-isotopic ribbon surfaces if and only if they are related by a finite sequence of moves $s_1$ to $s_{26}$ in Figures 3, 6–8.*

*Proof.* The 'if' part is trivial, since all the moves in Figures 6–8 represent special three-dimensional diagram isotopies. For the 'only if' part, we need to show that these moves do generate any three-dimensional diagram isotopy between planar diagrams.

Moves $s_5$ and $s_6$ allow us to restrict attention to special planar diagrams. Moreover, all the moves of Figures 7 and 8 contain only terminal ribbon intersections; hence, they can be performed in the context of special planar diagrams.

Now, consider two special planar diagrams representing ribbon surfaces $S_0$ and $S_1$, whose three-dimensional diagrams are isotopic in $R^3$, and let $H: R^3 \times [0,1] \to R^3$ be a smooth ambient isotopy such that $h_1(S_0) = S_1$.

For $i = 0, 1$, let $P_i$ be a simple spine of $S_i$, and $T_i = \pi(P_i)$ be the core of its diagram. Up to moves, we can assume that $h_1(T_0) = T_1$. Indeed, by cutting $S_1$ along the ribbon intersection arcs, we get an embedded surface $\widehat{S}_1 \subset R^3$ with some marked arcs, one in the interior and two along the boundary, for each ribbon intersection. This operation transforms the graphs $T_1$ and $h_1(T_0)$ into simple spines of $\widehat{S}_1$ relative to the marked arcs. Figure 9 shows the effect of the cut at the ribbon intersections in Figure 4.

From the intrinsic point of view, that is, considering $\widehat{S}_1$ as an abstract surface and forgetting its inclusion in $R^3$, the theory of simple spines implies that moves $s_7$, $s_8$ and the composition of the moves $s_5$ and $s_6$ suffice to transform $h_1(T_0)$ into $T_1$. In particular, the first two moves correspond to the well known moves for simple spines of surfaces, while the third together with the moves in Figure 5, which do not change the surface, relates the different positions of the spine with respect to the marked arcs in the interior of $\widehat{S}_1$. It remains only to observe that, up to the other moves in Figures 7 and 8, the portion of the surface involved in each single spine modification can be isolated in the planar diagram, as needed to perform the above-mentioned moves.

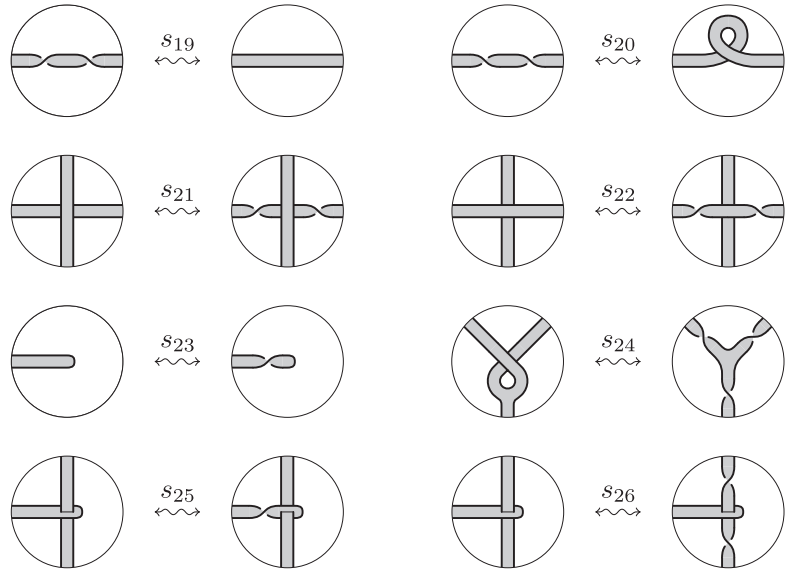

FIGURE 8. *Half-twisted isotopy moves.*

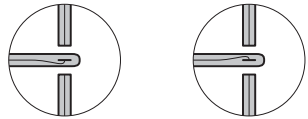

FIGURE 9. *Cutting the three-dimensional diagram at a ribbon intersection.*

So, let us suppose that $h_1(S_0, T_0) = (S_1, T_1)$. Note that the intermediate pairs $(S_t, T_t) = h_t(S_0, T_0)$ with $0 < t < 1$ do not necessarily project into special planar diagrams in $R^2$.

By transversality, we can assume that the graph $T_t$ regularly projects to a diagram in $R^2$ for every $t \in [0, 1]$, except a finite number of $t$'s corresponding to extended Reidemeister moves for graphs. For such exceptional $t$'s, the lines tangent to $T_t$ at its vertices are assumed not to be vertical.

We define $\Gamma \subset T_0 \times [0, 1]$ as the subspace of pairs $(x, t)$ for which the plane $T_{x_t} S_t$ tangent to $S_t$ at $x_t = h_t(x)$ is vertical (if $x \in T_0$ is a singular vertex, there are two such tangent planes and we require that one of them is vertical).

By a standard transversality argument, we can perturb $H$ in such a way that:

(a) $\Gamma$ is a graph embedded in $T_0 \times [0, 1]$ as a smooth stratified subspace of constant codimension 1 and the restriction $\eta : \Gamma \to [0, 1]$ of the height function $(x, t) \mapsto t$ is a Morse function on each edge of $\Gamma$;

(b) the edges of $\Gamma$ locally separate regions consisting of points $(x,t)$ for which the projection of $S_t$ into $R^2$ has opposite local orientations at $x_t$;

(c) the two planes tangent to any $S_t$ at a singular vertex of $T_t$ are not both vertical, and if one of them is vertical then it does not contain both the lines tangent to $T_t$ at that vertex.

As a consequence of (b), for each flat vertex $x \in T_0$ of valency one or three there are finitely many points $(x,t) \in \Gamma$, all of which have the same valency one or three as vertices of $\Gamma$. Similarly, as a consequence of (c), for each singular vertex $x \in T_0$ there are finitely many points $(x,t) \in \Gamma$, all of which have valency one or two as vertices of $\Gamma$. Moreover, the above-mentioned vertices of $\Gamma$ of valency one or three are the only vertices of $\Gamma$ of valency $\neq 2$.

Let $0 < t_1 < t_2 < \cdots < t_k < 1$ be the critical levels where one of the following holds:

(1) $T_{t_i}$ does not project regularly in $R^2$, because there is a point $x_i$ along an edge of $T_0$ such that the line tangent to $T_{t_i}$ at $h_{t_i}(x_i)$ is vertical;
(2) $T_{t_i}$ projects regularly in $R^2$, but its projection is not a graph diagram, due to a multiple tangency or crossing at some point;
(3) there is a point $(x_i, t_i) \in \Gamma$ with $x_i$ a univalent or a singular vertex of $T_0$;
(4) there is a critical point $(x_i, t_i)$ for the function $\eta$ along an edge of $\Gamma$.

Without loss of generality, we assume that only one of the four cases mentioned above occurs at any critical level $t_i$. Note that the points $(x,t)$ of $\Gamma$ such that $x \in T_0$ is a flat tri-valent vertex represent a subcase of case 2 and for this reason they are not included in case 3.

For $t \in [0,1] - \{t_1, t_2, \ldots, t_k\}$, there exists a sufficiently small regular neighbourhood $N_t$ of $T_t$ in $S_t$, such that the pair $(N_t, T_t)$ projects to a planar diagram.

We observe that the planar diagram of $N_t$ is uniquely determined up to diagram isotopy by (the diagram of) its core $T_t$ and by the tangent planes of $S_t$ at $T_t$. In fact, the half-twists of $N_t$ along the edges of $T_t$ correspond to the transversal intersections of $\Gamma$ with $T_0 \times \{t\}$ and their signs depend only on the local behaviour of the tangent planes of $T_t$. In particular, the planar diagrams of $(N_0, T_0)$ and $(N_1, T_1)$ coincide, up to planar diagram isotopy, with the original ones of $(S_0, T_0)$ and $(S_1, T_1)$.

If an interval $[t', t'']$ does not contain any critical level $t_i$, then each single half-twist persists between the levels $t'$ and $t''$, hence the planar isotopy relating the diagrams of $T_{t'}$ and $T_{t''}$ also relate the diagrams of $N_{t'}$ and $N_{t''}$, except for possible slidings of half-twists along ribbons over/under crossings. These can be realized by using moves $s_{19}$, $s_{21}$ and $s_{22}$.

At this point, the only thing left to show is that $N_{t'}$ and $N_{t''}$ are related by moves for $[t', t'']$ a sufficiently small neighbourhood of a critical level $t_i$. We do that separately for the four different types of critical levels.

If $t_i$ is of type 1, then a kink is appearing or disappearing along an edge of the core graph. When this kink is positive, the diagrams of $N_{t'}$ and $N_{t''}$ are directly related by move $s_{20}$ if $(x_i, t_i)$ is a local maximum point for $\eta$ and the kink is appearing, or if $(x_i, t_i)$ is a local minimum point for $\eta$ and the kink is disappearing. All the other cases of a positive kink can be reduced to the previous ones, by means of move $s_{19}$. The case of a negative kink is symmetric, we can just use moves $\bar{s}_i$ in place of the moves $s_i$.

If $t_i$ is of type 2, then either a regular isotopy move is occurring between $T_{t'}$ and $T_{t''}$ or two tangent lines at a tri-valent vertex $x_i$ of $T_{t_i}$ project to the same line in the plane. In the first case, the regular isotopy move occurring between $T_{t'}$ and $T_{t''}$, trivially extends to one of the moves $s_9$ to $s_{16}$. In the second case, $x_i$ may be either a flat or a singular vertex of $T_{t_i}$. If $x_i$ is a flat vertex, then the tangent plane to $S_t$ at $H(x_i, t)$ is vertical for $t = t_i$ and its projection reverses the orientation when $t$ passes from $t'$ to $t''$. Moves $s_{24}$ and $\bar{s}_{24}$ (modulo moves $s_9$ and $s_{19}$) describe the effect on the diagram of such a reversion of the tangent plane. If $x_i$ is a singular vertex, then $N_{t'}$ changes into $N_{t''}$ by one move $s_{17}$, $\bar{s}_{17}$, $s_{18}$ or $\bar{s}_{18}$.

If $t_i$ is of type 3, then either a half-twist is appearing/disappearing at the tip of the tongue of the surface corresponding to a univalent vertex or one of the two bands at the ribbon intersection corresponding to a singular vertex is being reversed in the plane projection. The first case corresponds to move $s_{23}$ or $\bar s_{23}$, while the second case corresponds, up to move $s_{19}$, to one of moves $s_{25}$, $\bar s_{25}$, $s_{26}$ or $\bar s_{26}$ (depending on the type of ribbon intersection and on which band is being reversed).

Finally, if $t_i$ is of type 4, a pair of cancelling half-twists is appearing or disappearing along a band, just as in move $s_{19}$. □

In the following we will focus on *flat planar diagrams*, meaning planar diagrams without half-twists. In other words, these are planar diagrams locally modelled on the spots (a), (b), (e) and (f) in Figure 2 (and possibly (g) and (h) in Figure 4).

Of course, only orientable ribbon surfaces can be represented by flat planar diagrams. In fact, a ribbon surface with a flat planar diagram has a preferred orientation induced by the projection in the plane of the diagram, which in this case is a regular map. Actually, any oriented ribbon surface is known to admit a flat planar diagram. But we will not need this fact here, and we just refer to [20] for its proof.

In contrast, finding a complete set of local moves representing 1-isotopy between oriented ribbon surfaces in terms of flat planar diagrams seems not to be so easy. These should include all the moves $s_1$ to $s_{18}$ in Figures 3, 6 and 7 and flat versions of some of the moves $s_{19}$ to $s_{26}$ in Figure 8.

However, this problem can be circumvented when using labelled orientable ribbon surfaces to represent branched coverings of $B^4$, thanks to the presence of the covering moves introduced in Section 5 (cf. Figure 24).

## 3. Braided surfaces

A regularly embedded smooth compact surface $S \subset B^2 \times B^2$ is called a (simply) *braided surface* of degree $m$ if the projection $\pi: B^2 \times B^2 \to B^2$ onto the first factor restricts to a simple branched covering $p = \pi_| : S \to B^2$ of degree $m$.

This means that there exists a finite set $A = \{a_1, a_2, \ldots, a_n\} \subset \operatorname{Int} B^2$ of branch points, such that the restriction $p_| : S - p^{-1}(A) \to B^2 - A$ is an ordinary covering of degree $m$, while over any branch point $a_i \in A$ there is only one singular point $s_i \in p^{-1}(a_i) \subset S$ and $p$ has local degree 2 at $s_i$, being locally smoothly equivalent to the complex map $z \mapsto z^2$.

For any singular point $s_i \in S$, there are local complex coordinates on the two factors centred at $s_i$, with respect to which $S$ has local equation $z_1 = z_2^2$. Actually, if we insist that those local coordinates preserve standard orientations, then we have two different possibilities, up to ambient isotopy, for the local equation of $S$ at $s_i$, namely $z_1 = z_2^2$ or $\bar z_1 = z_2^2$. We call $s_i$ a *positive twist point* for $S$ in the first case and a *negative twist point* for $S$ in the second case.

By a *braided isotopy* between the two braided surfaces $S, S' \subset B^2 \times B^2$, we mean a smooth ambient isotopy $(h_t)_{t \in [0,1]}$ of $B^2 \times B^2$ such that $h_1(S) = S'$ and each $h_t$ preserves the vertical fibres (those of the projection $\pi$); in other words, there exists a smooth ambient isotopy $(k_t)_{t \in [0,1]}$ of $B^2$ such that $\pi \circ h_t = k_t \circ \pi$ for every $t \in [0,1]$. In particular, if such a braided isotopy exists, then $S$ and $S'$ are isotopic through braided surfaces. Of course, braided isotopy reduces to *vertical isotopy* if $k_t = \operatorname{id}_{B^2}$ for every $t \in [0,1]$.

Now, assume that $* \in S^1$ is fixed once and for all as the base point of $B^2 - A$. Then, the classical theory of coverings tells us that the branched covering $p: S \to B^2$ is uniquely determined up to diffeomorphisms by the *monodromy* $\omega_p: \pi_1(B^2 - A) \to \Sigma_m$ of its restriction over $B^2 - A$ (defined only up to conjugation in $\Sigma_m$, depending on the numbering of the sheets).

Similarly, the braided surface $S \subset B^2 \times B^2$ is uniquely determined up to vertical isotopy by its *braid monodromy*, that is a suitable lifting $\omega_S: \pi_1(B^2 - A) \to \mathcal{B}_m$ of $\omega_p$ to the braid group $\mathcal{B}_m$ of degree $m$. This is defined in the following way: we take $\widetilde{*} = (*_1, *_2, \ldots, *_m) = p^{-1}(*) \subset \{*\} \times B^2 \cong B^2$ as the base point of the configuration space $\Gamma_m B^2$ of $m$ points in $B^2$, then for any $[\lambda] \in \pi_1(B^2 - A, *)$ we put $\omega_S([\lambda]) = [\widetilde{\lambda}] \in \pi_1(\Gamma_m B^2, \widetilde{*}) \cong \mathcal{B}_m$, where $\widetilde{\lambda}$ is the loop given by $\widetilde{\lambda}(t) = p^{-1}(\lambda(t)) \subset \{\lambda(t)\} \times B^2 \cong B^2$ for any $t \in [0, 1]$. We can immediately see that $\sigma \circ \omega_S = \omega_p$, where $\sigma: \mathcal{B}_m \to \Sigma_m$ is the canonical homomorphism giving the permutation associated to a braid. Like the monodromy $\omega_p$, the braid monodromy $\omega_S$ is defined only up to conjugation in $\mathcal{B}_m$, depending on the identification $\pi_1(\Gamma_m B^2, \widetilde{*}) \cong \mathcal{B}_m$.

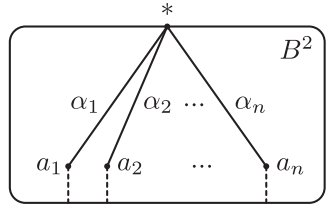

FIGURE 10. *The standard Hurwitz system.*

The local model of the twists points forces the braid monodromy $\omega_S(\mu)$ of any meridian $\mu \in \pi_1(B^2 - A)$ around a branch point $a \in A$ to be a half-twist $\beta^{\pm 1} \in \mathcal{B}_n$ around an arc $b \subset B^2$ between two points $*_j$ and $*_k$. The arc $b$ turns out to be uniquely determined up to ambient isotopy of $B^2$ mod $\widetilde{*}$, while the half-twist $\beta^{\pm 1}$ is positive (right-handed) or negative (left-handed) according to the sign of the twist point $s_i \in S$.

Conversely, as we will see shortly, any homomorphism $\varphi: \pi_1(B^2 - A) \to \mathcal{B}_m$ that sends meridians around the points of $A$ to positive or negative half-twists around intervals in $B^2$ is the braid monodromy of a braided surface $S \subset B^2 \times B^2$ with branch set $A$.

We recall that $\pi_1(B^2 - A)$ is freely generated by any set of meridians $\alpha_1, \alpha_2, \ldots, \alpha_n$ around the points $a_1, a_2, \ldots, a_n$, respectively. An ordered sequence $(\alpha_1, \alpha_2, \ldots, \alpha_n)$ of such meridians is called a *Hurwitz system* for $A$ when the following properties hold: (i) each $\alpha_i$ is realized as a counterclockwise parametrization of the boundary of a regular neighbourhood in $B^2$ of a non-singular arc from $*$ to $a_i$, which we still denote by $\alpha_i$; (ii) except for their end points, the arcs $\alpha_1, \alpha_2, \ldots, \alpha_n$ are pairwise disjoint and contained in Int $B^2 - A$; (iii) around $*$ the arcs $\alpha_1, \alpha_2, \ldots, \alpha_n$ appear in the counterclockwise order, so that the composition loop $\alpha_1 \alpha_2 \cdots \alpha_n$ is homotopic in $B^2 - A$ to the usual counterclockwise generator $\alpha \in \pi_1(S^1)$, and the points of $A$ are assumed to be indexed accordingly. Up to ambient isotopy of $B^2$ fixing $S^1$ but not $A$, any Hurwitz system looks like the standard one depicted in Figure 10.

For the sake of convenience, here the disc $B^2$ is drawn as $B^1 \times B^1$ with rounded corners. *Actually, in all the pictures, we will always draw both the horizontal and the vertical fibres of $B^2 \times B^2$ as $B^1 \times B^1$ with rounded corners.*

There is a natural transitive action of the braid group $\mathcal{B}_n \cong \pi_1(\Gamma_n \operatorname{Int} B^2, A)$ on the set of Hurwitz systems for $A$. To any such Hurwitz system $(\alpha_1, \alpha_2, \ldots, \alpha_n)$ we associate a set of standard generators $\xi_1, \xi_2, \ldots, \xi_{n-1}$ of $\mathcal{B}_n$, with $\xi_i$ the right-handed half-twists around the interval $x_i \simeq \bar{\alpha}_i \alpha_{i+1}$ with end points $a_i$ and $a_{i+1}$. Under the action of $\mathcal{B}_n$, each $\xi_i$ transforms $(\alpha_1, \alpha_2, \ldots, \alpha_n)$ into $(\alpha'_1, \alpha'_2, \ldots, \alpha'_n)$ with $\alpha'_i = \alpha_i \alpha_{i+1} \alpha_i^{-1}$, $\alpha'_{i+1} = \alpha_i$ and $\alpha'_k = \alpha_k$ for $k \neq i, i+1$. This will be referred to as the *ith elementary transformation* $\xi_i$ (cf. Figure 11). It turns out that any two Hurwitz systems for $A$ are related by a finite number of consecutive elementary transformations $\xi_i^{\pm 1}$ with $i = 1, 2, \ldots, n-1$.

Given a Hurwitz system $(\alpha_1, \alpha_2, \ldots, \alpha_n)$ for $A$, we can represent the braid monodromy of the braided surface $S$ by the sequence $(\beta_1 = \omega_S(\alpha_1), \beta_2 = \omega_S(\alpha_2), \ldots, \beta_n = \omega_S(\alpha_n))$ of positive or negative half-twists in $\mathcal{B}_m$, and the monodromy of the branched covering $p: S \to B^2$ by the sequence $(\tau_1 = \sigma(\beta_1), \tau_2 = \sigma(\beta_2), \ldots, \tau_n = \sigma(\beta_n))$ of the associated transpositions in $\Sigma_m$.

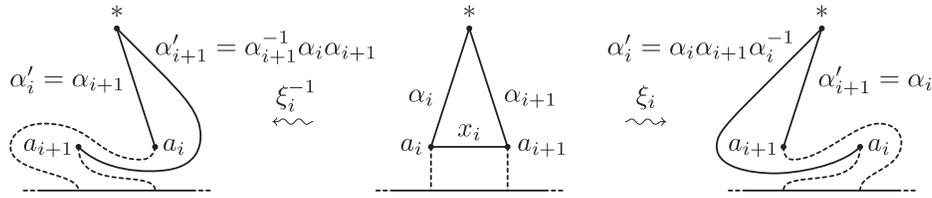

FIGURE 11. *Elementary transformations.*

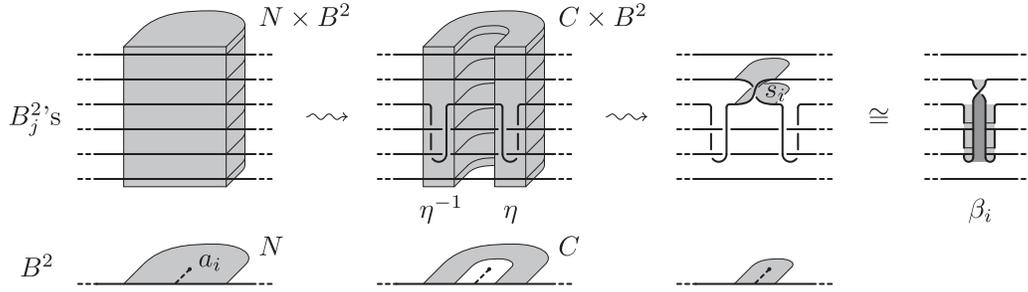

FIGURE 12. *From braid monodromy to braided surfaces.*

Conversely, starting from any sequence $(\beta_1, \beta_2, \ldots, \beta_n)$ of positive or negative half-twists in $\mathcal{B}_m$, we can construct a braided surface $S = S(m; \beta_1, \beta_2, \ldots, \beta_n)$ of degree $m$, whose braid monodromy is determined by $\beta_i = \omega_S(\alpha_i)$, as follows (cf. [20, Section 2]). First, we fix a base point $\widetilde{*} = (*_1, *_2, \ldots, *_m) \in \Gamma_m B^2$ and consider the $m$ horizontal copies of $B^2$ given by $B_j^2 = B^2 \times \{*_j\} \subset B^2 \times B^2$ for any $j = 1, 2, \ldots, m$. We assume that the points $*_j$ form an increasing sequence in $B^1 \subset B^2$, to let the discs $B_1^2, B_2^2, \ldots, B_m^2$ appear to be stacked up on the top of each other in that order, when we look at the three-dimensional picture given by the canonical projection $\pi : B^2 \times B^2 \to B^2 \times B^1$. Then, we consider the standard generators $\xi_1, \xi_2, \ldots, \xi_{m-1}$ of $\mathcal{B}_m \cong \pi_1(\Gamma_m B^2, \widetilde{*})$, with $\xi_i$ the right-handed half-twist around the vertical interval $x_i$ between $*_i$ and $*_{i+1}$. Finally, we choose a family $\delta_1, \delta_2, \ldots, \delta_n$ of disjoint arcs in $B^2$, respectively, joining the points $a_1, a_2, \ldots, a_n$ to $S^1$, like the dashed ones in Figure 10, which form a splitting complex for the branched covering $p : S \to B^2$, and we do the following for each $i = 1, 2, \ldots, n$: (i) we express the half-twist $\beta_i$ as $\eta^{-1} \xi_{j_i}^{\pm 1} \eta$, with $\xi_{j_i}$ a standard generator of $\mathcal{B}_m$ and $\eta \in \mathcal{B}_m$ such that $\eta(x_{j_i}) = b_i$ is the arc around which $\beta_i$ is defined, thanks to the transitive action of $\mathcal{B}_m$ on the set of arcs between points of $A$; (ii) we deform the discs $B_j^2$ by a vertical ambient isotopy supported inside $N \times B^2$ for a small regular neighbourhood $N$ of $\delta_i$ in $B^2$, which realizes the braid $\eta$ over each fibre of a collar $C \subset N$ of the boundary of $N$ in $B^2$, while it does not depend on the first component over $N - C$; (iii) we replace the two adjacent discs $(N - C) \times \{*_{j_i}\} \subset B_{j_i}^2$ and $(N - C) \times \{*_{j_i+1}\} \subset B_{j_i+1}^2$ by the local model described above for a positive or negative twist point, depending on the sign of the half-twist $\beta_i$ (that is on the exponent of $\xi_{j_i}^{\pm 1}$).

Up to horizontal isotopy, the last construction results in attaching to the horizontal discs a narrow half-twisted vertical band, which we still denote by $\beta_i$, as the relative half-twist of the starting sequence. The band $\beta_i$ has a half-twist whose sign is the opposite of that of the original half-twist $\beta_i$ (and of the twist-point $s_i \in S$), hence it contributes to the boundary braid of the three-dimensional picture a half-twist having the same sign as the original one. Moreover, the core of the band $\beta_i$ is the arc $b_i$ around which the half-twist $\beta_i$ was defined, translated to the fibre over $a_i$. See Figure 12 for the case when $\beta_i$ is the negative half-twist $\eta^{-1} \xi_4^{-1} \eta \in \mathcal{B}_6$ around the arc $b_i = \eta(x_4)$ with $\eta = \xi_3 \xi_2^2 \xi_3$.

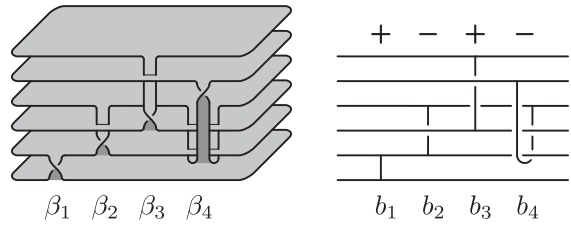

Figure 13. *Band presentation and line diagram of a braided surface.*

The identification $B^2 \times B^2 - \{(0, *)\} \cong B^4 - \{\infty\} \cong R_+^4$ given by a suitable rounding (smoothing the corners) of $B^2 \times B^2$ followed by the standard orientation preserving conformal equivalence $B^4 \cong R_+^4 \cup \{\infty\}$, makes the braided surface $S$ we have just constructed into a ribbon surface $\widehat{S} \subset R_+^4$. In fact, the projection of $\widehat{S}$ in $R^3$ turns out to be a three-dimensional diagram provided that the images of the arcs $\widehat{b}_i$ meet transversally those of the discs $\widehat{B}_j^2$, each ribbon intersection arc being formed by a band $\widehat{\beta}_i$ passing through a disc $\widehat{B}_j^2$, in correspondence with a transversal intersection point between $\widehat{b}_i$ and $\widehat{B}_j^2$. Therefore, the 1-handlebody decomposition of $\widehat{S}$ given by the discs $\widehat{B}_1^2, \widehat{B}_2^2, \ldots, \widehat{B}_m^2$ (as the 0-handles) and by the bands $\widehat{\beta}_1, \widehat{\beta}_2, \ldots, \widehat{\beta}_n$ (as the 1-handles), turns out to be an adapted one.

We call the ribbon surface $\widehat{S} \subset R_+^4$ a *band presentation* of the braided surface $S \subset B^2 \times B^2$. For example, on the left side of Figure 13 we see a band presentation of the braided surface $S(6; \beta_1, \beta_2, \beta_3, \beta_4)$ of degree 6 arising from the sequence $(\beta_1 = \xi_1, \beta_2 = \xi_3^{-1}\xi_2^{-1}\xi_3, \beta_3 = \xi_5^{-1}\xi_4\xi_3\xi_4^{-1}\xi_5, \beta_4 = \xi_3^{-1}\xi_2^{-2}\xi_3^{-1}\xi_4^{-1}\xi_3\xi_2^2\xi_3)$ of half-twists in $\mathcal{B}_6$.

A more economical way to represent braided surfaces in terms of band presentations is provided by *line diagrams*. These are just a variation in the charged fence diagrams introduced by Rudolph [22]. Namely, they consist of $m$ horizontal lines standing for the discs $\widehat{B}_1^2, \widehat{B}_2^2, \ldots, \widehat{B}_m^2$ and $n$ arcs between them given by the cores $\widehat{b}_1, \widehat{b}_2, \ldots, \widehat{b}_n$ of the bands $\widehat{\beta}_1, \widehat{\beta}_2, \ldots, \widehat{\beta}_n$, with the signs of the corresponding half-twists on the top. The right side of Figure 13 shows the line diagram of the surface depicted on the left side.

Of course, a braided surface $S$ has different band presentations, depending on the choice of various objects involved in the construction above: (i) the Hurwitz system $(\alpha_1, \alpha_2, \ldots, \alpha_n)$; (ii) the base point $\widetilde{*} \in \Gamma_m B^2$; (iii) the particular realizations of the arcs $b_1, b_2, \ldots, b_n$ within their isotopy classes.

The choices at points 2 and 3 are not relevant up to vertical isotopy of the braided surface $S(m; \beta_1, \beta_2, \ldots, \beta_m)$, but still they can affect the ribbon surface diagram of the corresponding band presentation.

Concerning point 1, we observe that different Hurwitz systems lead to different sequences of half-twists. As any two Hurwitz systems are related by elementary transformations and their inverses, the same holds for the corresponding sequences of half-twists.

Adopting Rudolph's terminology [20], we call such an elementary transformation of the sequence of half-twists a *band sliding*. Namely, the sliding of $\beta_{i+1}$ over $\beta_i$ changes the sequence $(\beta_1, \beta_2, \ldots, \beta_n)$ into $(\beta_1', \beta_2', \ldots, \beta_n')$, with $\beta_i' = \beta_i\beta_{i+1}\beta_i^{-1}$, $\beta_{i+1}' = \beta_i$ and $\beta_k' = \beta_k$ for $k \neq i, i+1$. The inverse transformation is the sliding of $\beta_i'$ over $\beta_{i+1}'$. Actually, these can be geometrically interpreted as genuine embedded 1-handle slidings only in the case when $b_i$ and $b_{i+1}$ can be realized as arcs whose intersection is one of their end points. On the other hand, it reduces to the interchange of $\beta_i$ and $\beta_{i+1}$ if $b_i$ and $b_{i+1}$ can be realized as disjoint arcs, hence the two half-twists commute. Figure 14 shows a band interchange followed by a geometric band sliding, in terms of band presentations and line diagrams.

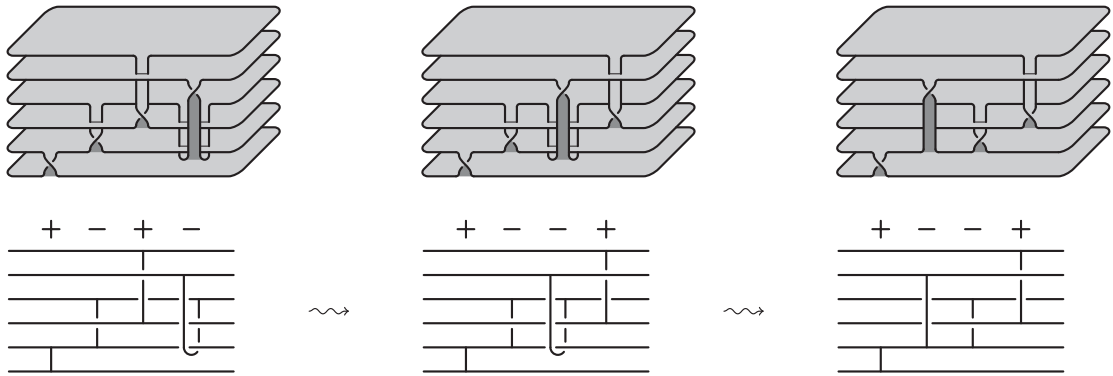

Figure 14. *Band interchange and sliding.*

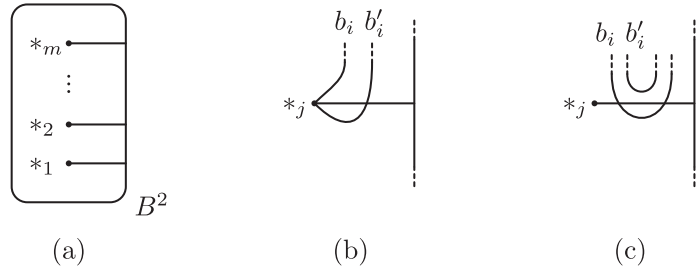

Figure 15. *Isotoping an arc $b_i$.*

Recalling that any two Hurwitz systems for a given branch set $A$ are isotopically equivalent to the standard one (if we do not insist on keeping $A$ fixed), when considering braided surfaces up to braided isotopy we can always assume the Hurwitz system to be the standard one. From this point of view, we can say that a sequence of half-twists $(\beta_1, \beta_2, \ldots, \beta_n)$ in $\mathcal{B}_m$, without any reference to a specific Hurwitz system, uniquely determines the braided surface $S(m; \beta_1, \beta_2, \ldots, \beta_n)$ up to braided isotopy. Moreover, the braided surfaces determined by two such monodromy sequences are braided isotopic if and only if they are related by simultaneous conjugation of all the $\beta_i$'s in $\mathcal{B}_m$ and band slidings (hence cyclic shift of the $\beta_i$'s as well).

PROPOSITION 2. *All the band presentations of a braided surface $S$ are 1-isotopic. Moreover, if $S'$ is another braided surface related to $S$ by a braided isotopy, then the band presentations of $S'$ are 1-isotopic to those of $S$.*

*Proof.* We first address the dependence of the band presentation of $S$ on the arcs $b_1, b_2, \ldots, b_n$, assuming that the Hurwitz system is fixed. It is clear from the construction of $\widehat{S}$ that the ribbon intersections of the band $\widehat{\beta}_i$ with the horizontal discs $\widehat{B}_j^2$ arise from the (transversal) intersections of $b_i$ with the horizontal arcs joining the points $*_j$ with $\mathrm{Bd}\, B^2$ depicted in Figure 15(a). On the other hand, we recall that $b_i$ is uniquely determined up to ambient isotopy of $B^2$ mod $\widetilde{*}$. By transversality, we can assume that such an isotopy essentially modifies the intersections of $b_i$ with those horizontal arcs only at a finite number of levels, when $b_i$ changes to $b'_i$ as in Figure 15(b) or (c) up to symmetry. Then, except for these critical levels any isotopy of the arc $b_i$ induces a three-dimensional diagram isotopy of $\widehat{S}$, while the modifications induced on $\widehat{S}$ at the critical levels of type (b) and (c) can be realized by straightforward applications of the 1-isotopy moves $s_1$ and $s_{2,3}$, respectively.

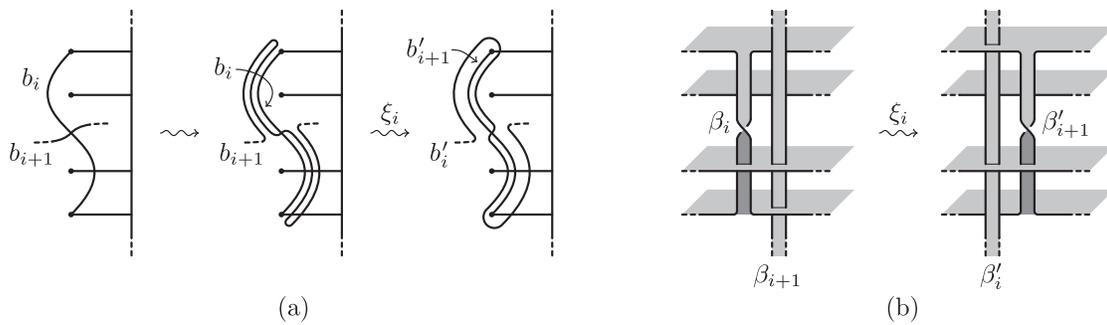

FIGURE 16. *Sliding $\beta_{i+1}$ over a positive $\beta_i$ (the non-trivial case).*

At this point, having proved the independence on the arcs $b_i$, we observe that the vertical isotopy relating the braided surfaces resulting from different choices of the base point $\widetilde{*} \in \Gamma_m B^2$ (subject to the condition the $*_i$'s form an increasing sequence in $B^1 \subset B^2$) induces three-dimensional diagram isotopy on the band presentation.

For the dependence of the band presentation on the Hurwitz system and for the second part of the proposition, it suffices to consider the case of the elementary transformation of the sequence of half-twists $(\beta_1, \beta_2, \ldots, \beta_n)$ given by the sliding of $\beta_{i+1}$ over $\beta_i$. If the arcs $b_i$ and $b_{i+1}$ are disjoint, hence $\beta_i$ and $\beta_{i+1}$ commute, the band presentation only changes by a three-dimensional diagram isotopy. In the case when $b_i \cap b_{i+1}$ consists of one common end point, we have a true embedded sliding, which can be easily realized by the 1-isotopy moves $s_{2,3}$. The case when $b_i$ and $b_{i+1}$ share both end points and nothing else is similar, being reducible to two consecutive true embedded slidings. Then, we are left to consider the case when $b_i$ and $b_{i+1}$ have some transversal intersection point (possibly in addition to some common end point). In this case, we first isotope the arc $b_{i+1}$ so that each transversal intersection is contained in a portion of the arc that runs nearly parallel to all the arc $b_i$. There are essentially two different ways to do that, the right one depending on the sign of the half-twist $\beta_i$. The first step of Figure 16(a) shows how to deal with a single transversal intersection for a positive $\beta_i$ ($\beta_{i+1}$ should be isotoped in the other way for a negative $\beta_i$). In any case, according to the first part of the proof, isotoping $\beta_{i+1}$ induces 1-isotopy on the band presentation $\widehat{S}$. After that, the desired elementary transformation amounts to passing the band $\beta_{i+1}$ through the band $\beta_i$ in the band presentation $\widehat{S}$, as it can be easily realized by looking again at the example described in Figure 16 (in (b) only the portion of $\beta_{i+1}$ parallel to $\beta_i$ is shown). Then, to conclude the proof, it suffices to notice that $\beta_{i+1}$ can be passed through $\beta_i$ by means of a sequence of 1-isotopy moves $s_{2,3,4}$. In particular, move $s_4$ is needed to pass the ribbon intersections of $\beta_i$ with the discs $B_j^2$. □

In the following, we will not distinguish between a braided surface $S$ and any band presentation of it, taking into account that there is a canonical identification between them and that the latter is uniquely determined up to 1-isotopy.

A braided surface $S$ of degree $m$ with $n$ twist points can be deformed to a braided surface $S'$ of degree $m+1$ with $n+1$ twist points, called an *elementary stabilization* of $S$, by expanding a new half-twisted band from one of the horizontal discs of $S$ and then a new horizontal disc from the tip of that band. Looking at the three-dimensional diagram, we see that the band presentations of $S$ and $S'$ are 1-isotopic. In fact, apart from the three-dimensional diagram isotopy only move $s_2$ is needed when the new band is pushed through a disc. The inverse process, that is cancelling a band $\beta_i$ and a horizontal disc $B_j^2$ from a braided surface $S$ to get an *elementary destabilization* of it, can be performed when $\beta_i$ is the only band attached to $B_j^2$ and no band is linked with (passes through, in the three-dimensional diagram) $B_j^2$. For example, in the braided surface of Figure 13 the band $\beta_3$ can be cancelled with the disc $B_6^2$, and after that (but not before) the band $\beta_4$ can be cancelled with the disc $B_5^2$. A *(de)stabilization* is the result of consecutive elementary (de)stabilizations.

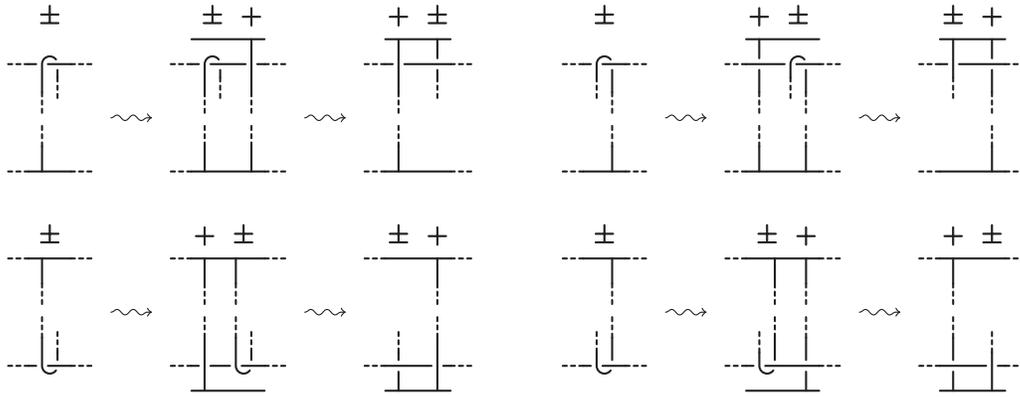

Figure 17. *Making bands monotonic.*

We say that a band $\beta_i$ of a braided surface $S$ is a *monotonic band*, if it has the form $\xi_{k-1}^{-\varepsilon_{k-1}} \xi_{k-2}^{-\varepsilon_{k-2}} \cdots \xi_{j+1}^{-\varepsilon_{j+1}} \xi_j^{\pm 1} \xi_{j+1}^{\varepsilon_{j+1}} \cdots \xi_{k-2}^{\varepsilon_{k-2}} \xi_{k-1}^{\varepsilon_{k-1}}$ for some $j < k$ and $\varepsilon_h = \pm 1$. In other words, $\beta_i$ appears to run monotonically (with respect to the coordinate $x_3$) from $B_j^2$ to $B_k^2$ in the three-dimensional diagram of $S$, and its core $b_i$ can be drawn as a vertical segment in the line diagram of $S$ (remember that we are assuming the standard generators $\xi_1, \xi_2, \ldots, \xi_{m-1}$ of $\mathcal{B}_m$ to be half-twists around vertical arcs). For example, the bands $\beta_1, \beta_2$ and $\beta_3$ in Figure 13 are monotonic, whereas $\beta_4$ is not. $S$ is called a braided surface *with monotonic bands* if all its bands are monotonic. In the following proposition, we see that stabilization and band sliding enable us to transform any braided surface into one with monotonic bands.

PROPOSITION 3. *Any braided surface $S$ admits a positive stabilization $S'$ with monotonic bands up to braided isotopy. Moreover, since stabilization is realizable by 1-isotopy, $S$ and $S'$ are 1-isotopic.*

*Proof.* In Figure 17 we see how to eliminate the first extremal point along the core $b$ of a band $\beta$ (the one having sign $\pm$ in the diagrams) by a suitable positive elementary stabilization and the subsequent sliding of the band $\beta$ over the new stabilizing band. For the sake of clarity, here all the four possible cases are shown, even if they are symmetric to each other. In all the cases, the new stabilizing band is a monotonic band that runs parallel to the first monotonic portion of $\beta$ (in particular, it passes through the same horizontal discs).

Iterating this process for all the extremal points along $b$, we can replace the band $\beta$ with a sequence of monotonic bands. Once this is done for all the bands of $S$, we get a braided surface with monotonic bands.

It remains to observe that all the elementary stabilizations can be performed at the beginning to obtain the desired stabilization $S'$, while leaving all the band slidings at the end to give a braided isotopy from $S'$ to a braided surface with monotonic bands. □

To conclude this section we observe that any braided surface $S$ is orientable, carrying the preferred orientation induced by the branched covering $p \colon S \to B^2$. Therefore, any band presentation of it admits a flat planar diagram. For a braided surface $S$ with monotonic bands, this can be easily obtained through the three-dimensional diagram isotopy given by the following simple procedure (cf. [**20**] and see Figure 18 for an example): first flatten the half-twisted bands by inserting a half-curl at their bottom ends, then contract the discs to non-overlapping horizontal bands. We call this the *flattening procedure*.

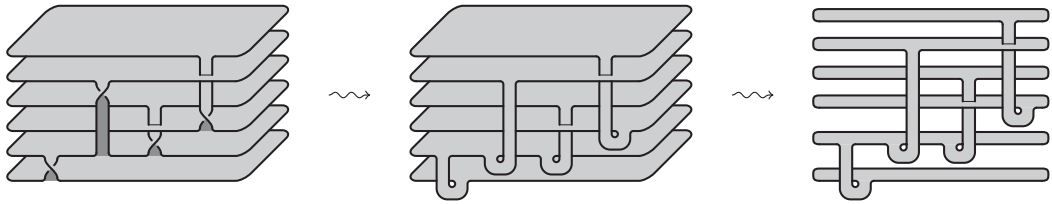

Figure 18. *Getting a flat planar diagram of a band presentation.*

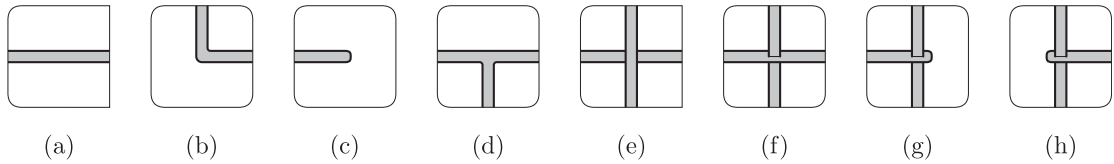

Figure 19. *Local models for rectangular diagrams.*

Conversely, Rudolph provided in [**20**] a *braiding procedure* to produce a three-dimensional diagram isotopy, which makes an orientable ribbon surface given by a flat planar diagram into a band presentation of a braided surface. In the next section, we will describe this braiding procedure in a revised form suitable for our purposes.

## 4. Rudolph's braiding procedure

Following [**20**], we start from the observation that up to planar ambient isotopy any flat planar diagram can be assumed to have all the bands parallel to the coordinate axes. The flat planar diagrams with this property will be the input for the braiding procedure. Before going on, let us give a more precise definition of them.

A *rectangular diagram* of a ribbon surface is a flat planar diagram, whose local configurations are those described in Figure 19, possibly rotated by $\pi/2$, $\pi$ or $3\pi/2$ radians. We denote by prime, double prime and triple prime, respectively, the configurations obtained by these rotations. In particular, (g) and (h) should be thought as contractions of (f) and (f$''$) juxtaposed with (c) and (c$''$), respectively. Arbitrarily many (possibly rotated) configurations of types (d), (e) and (f) can occur along any horizontal or vertical band, and (possibly rotated) configurations of types (b), (d), (g) and (h) can appear at both the ends of the band, but different horizontal bands are always assumed to have different ordinates, and different vertical bands are always assumed to have different abscissas.

Rectangular diagrams will always be considered up to plane ambient isotopy through diffeomorphisms of the form $(x, y) \mapsto (h_1(x), h_2(y))$ with $h_1$ and $h_2$ monotonic increasing real functions.

The reader may have noticed that in Figure 19 some of the corners of the boxes are rounded and some are not. We use this detail to specify the rotations we want to consider and admit, according to the following rule: a box can be rotated only in the positions such that the bottom-left corner is rounded. Of course, due to the symmetry of (a) and (e), this constraint is not effective here, but it will be in the next figures.

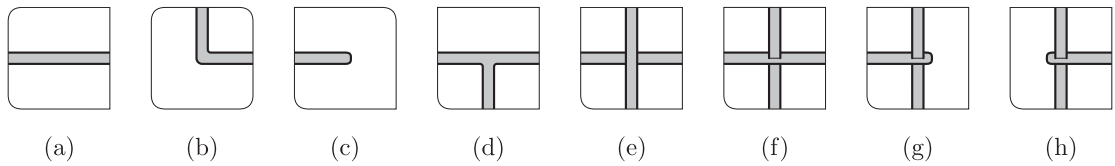

Figure 20. *Allowed local models for the restricted braiding procedure.*

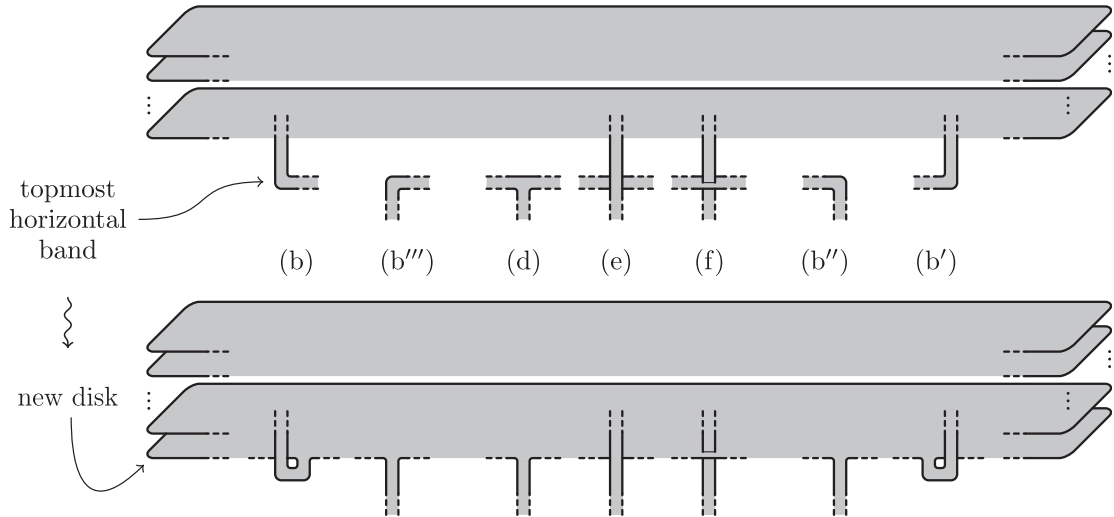

Figure 21. *The restricted braiding procedure (step 1).*

We first define a restricted version of the braiding procedure on the rectangular diagrams whose local configurations are constrained as in Figure 20, according to the above rule. Namely, the allowed local configurations are (a), (a′), (b), (b′), (b″), (b‴), (c), (c″), (d), (e), (f), (g) and (h). Also in this context (g) and (h) are just notational contractions and we will not consider them as separate cases. A rectangular diagram presenting only these local configurations is said to be in *restricted form*.

Starting from a rectangular diagram in restricted form, the first step of the braiding procedure is to transform each horizontal band, in the order from top to bottom, into a disc inserted under the previous ones. For a horizontal band, the left end may be of type (b), (b‴) or (c″), and in the first two cases we get a vertical band attached to the new disc as shown in Figure 21, while in the third case we do not get any vertical band. Analogously, the right end may be of type (b′), (b″) or (c), and in the first two cases we get a vertical band attached to the new disc as shown in Figure 21, while in the third case we do not get any vertical band. On the other hand, we have arbitrarily many bands attached to the new disc in correspondence with the local configurations like (d), and arbitrarily many bands passing either in front or through the new disc in correspondence the local configurations like (e) and (f), respectively, as shown in Figure 21.

After each horizontal band has been transformed into a disc, we isotope all the resulting half-curls to half-twists according to Figure 22. The final result is a band presentation of a braided surface with monotonic bands.

Now, to extend the braiding procedure to any rectangular diagram, we first replace all the local configurations not included in Figure 20, by means of the plane diagram isotopies described in Figure 23. This produces a rectangular diagram in restricted form, which depends on the order of the replacements. However, the diagrams obtained following different orders can be easily proved to be equivalent up to the moves $r_1$ and $r_1'$ defined in Section 8 (see Figure 39). Since this fact will suffice for our purposes, we do not need to worry about the order of replacements. For the moment, let us assume that these are performed in the lexicographic order from top to bottom and then from left to right.

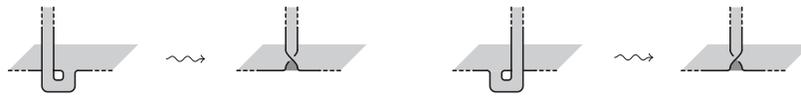

Figure 22. *The restricted braiding procedure (step 2).*

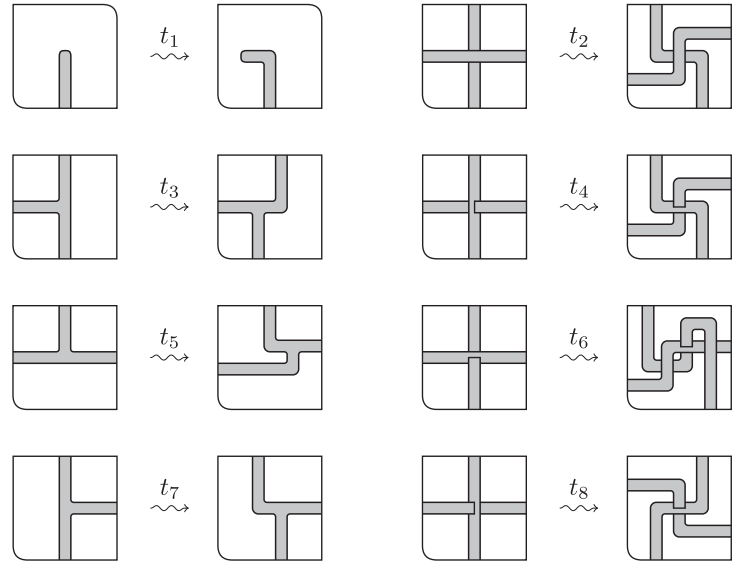

Figure 23. *Extending the braiding procedure to any rectangular diagram.*

Proposition 4. *The braiding procedure described above produces a three-dimensional diagram isotopy from any rectangular diagram of an orientable ribbon surface to a band presentation of a braided surface with monotonic bands. Moreover, any such band presentation can be obtained in this way, starting from the rectangular diagram given by the flattening procedure (defined at the end of the previous section) applied to it.*

*Proof.* The first part of the statement is clear from the construction above. In particular, each vertical band is created starting from the top when a configuration of type (d), (b″) or (b‴) is meet, then it keeps going down, possibly passing through some discs in correspondence with the configurations of type (f), until it ends at the bottom with a half-curl deriving from a configuration of type (b) or (b′). Then, the result of the braiding procedure is a braided surface with monotonic bands.

For the second part of the statement, it suffices to observe that the flattening procedure applied to a band presentation of a braided surface with monotonic bands produces a rectangular diagram in restricted form. Then, the braiding procedure applied to such diagram can be easily seen to give back the original band presentation. □

## 5. *Four-Dimensional 2-handlebodies*

By a *four-dimensional* 2-*handlebody* we mean a compact orientable four-manifold $W$ endowed with a handlebody structure, whose handles have indices at most 2. We call *2-equivalence* the equivalence relation on four-dimensional 2-handlebodies generated by 2-*deformations*, meaning handle isotopy, handle sliding and addition/deletion of cancelling pairs of handles of indices less than or equal to 2. Of course, 2-equivalent, four-dimensional 2-handlebodies are diffeomorphic, while the converse is not known and likely false.

Here, we consider four-dimensional 2-handlebodies as simple covers of $B^4$ branched over ribbon surfaces. We recall that a smooth map $p: W \to B^4$ is called a $d$-fold *branched covering* if there exists a smooth two-dimensional subcomplex $S \subset B^4$, the *branch set*, such that the restriction $p_|: W - p^{-1}(S) \to B^4 - S$ is a $d$-fold ordinary covering. We will always assume that $S$ is a ribbon surface in $R_+^4 \subset R_+^4 \cup \{\infty\} \cong B^4$. In this case, $p$ can be completely described in terms of the monodromy $\omega_p: \pi_1(B^4 - S) \to \Sigma_d$, by labelling each region of the three-dimensional diagram of $S$ with the permutation $\omega_p(\mu)$ associated to a meridian $\mu$ around it, in such a way that the usual Wirtinger relations at the crossings are respected. Conversely, any $\Sigma_d$-labelling of $S$ respecting such relations actually describes a covering of $B^4$ branched over $S$. Moreover, $p$ is called a *simple* branched covering if over any branch point $y \in S$ there is only one singular point $x \in p^{-1}(y)$ and $p$ has local degree 2 at $x$, being locally smoothly equivalent to the complex map $(z_1, z_2) \mapsto (z_1, z_2^2)$. In terms of the corresponding labelling, this means that each region is labelled by a transposition in $\Sigma_d$. We will refer to a ribbon surface with such a labelling by transpositions in $\Sigma_d$ as a *labelled ribbon surface*.

PROPOSITION 5. *A simple covering $p: W \to B^4$ branched over a ribbon surface determines a four-dimensional 2-handlebody decomposition $H_p$ of $W$, well defined up to 2-deformations.*

*Proof.* Following [**16**], once an adapted 1-handlebody decomposition $S = (D_1 \sqcup \cdots \sqcup D_m) \cup (B_1 \sqcup \cdots \sqcup B_n)$ of $S$ is given, with discs $D_i$ as 0-handles and bands $B_j$ as 1-handles, a 2-handlebody decomposition $W = (H_1^0 \sqcup \cdots \sqcup H_d^0) \cup (H_1^1 \sqcup \cdots \sqcup H_m^1) \cup (H_1^2 \sqcup \cdots \sqcup H_n^2)$, where $d$ is the degree of $p$, can be constructed as follows. We put $S_0 = D_1 \sqcup \cdots \sqcup D_m \subset B^4$ and denote by $p_0: W_1 \to B^4$ the $d$-fold simple covering branched over $S_0$ with the labelling inherited by $S$. Then, we put $W_1 = (H_1^0 \sqcup \cdots \sqcup H_d^0) \cup (H_1^1 \sqcup \cdots \sqcup H_m^1)$, where the 0-handles $H_1^0, \ldots, H_d^0 \cong B^4$ are the sheets of the covering $p_0$ and we have a 1-handle $H_i^1$ between the 0-handles $H_k^0$ and $H_l^0$ for each disc $D_i \subset S_0$ with label $(k\,l)$. Finally, $W$ can be obtained by attaching to $W_1$ a 2-handle $H_j^2$ for each band $B_i \subset S$, whose attaching map is described by the framed knot given by the unique annular component of $p_0^{-1}(B_j) \subset \mathrm{Bd}\, W_1$ (here we think $B_j \subset R^3$ as a band in the three-dimensional diagram of $S$). A detailed discussion of this construction in terms of Kirby diagrams can be found in [**3**, Section 2].

Now, according to [**3**, Proposition 2.2], the 2-equivalence class of the 2-handlebody decomposition of $W$ we have just described does not depend on the particular choice of the 1-handlebody decomposition of $S$. □

In light of the above proposition, it makes sense to say that any $\Sigma_d$-labelled ribbon surface $S \subset B^4$ representing a simple branched covering $p$, also represents the four-dimensional 2-handlebody $H_p$ up to 2-deformations.

In terms of this representation, the addition of a pair of cancelling 0- and 1-handles to the handlebody structure of $W$ can be interpreted as the addition of a $(d+1)$th extra sheet to the covering and the corresponding addition to $S$ of a separate trivial disc $D_{m+1}$ labelled $(i\,d+1)$ with $i \leqslant d$. We call *elementary stabilization* this operation, which changes a $d$-fold branched covering into a $(d+1)$-fold one representing the same handlebody up to 2-deformation, and

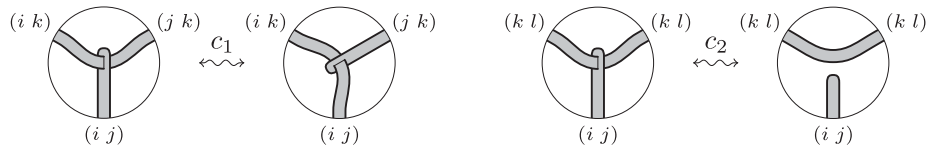

Figure 24. *The covering moves.*

*elementary destabilization* its inverse. Also in this context, by a *(de)stabilization* we mean the result of consecutive elementary (de)stabilizations.

On the other hand, the addition/deletion of a pair of cancelling 1- and 2-handles in the handlebody structure of $W$ can be interpreted as the addition/deletion of a corresponding cancelling disc and band in the handlebody structure of $S$. This leaves essentially unchanged the labelled ribbon surface $S$ (possibly up to some 1-isotopy moves $s_2$ occurring when the band passes through some discs), hence the covering $p: W \to B^4$ as well.

The following proposition summarizes results from [**3**, **16**].

PROPOSITION 6. *Up to 2-deformations, any connected four-dimensional 2-handlebody can be represented as a simple 3-fold branched covering of $B^4$, by a $\Sigma_3$-labelled ribbon surface in $B^4$. Two labelled ribbon surfaces in $B^4$ represent 2-equivalent connected four-dimensional 2-handlebodies if and only if, after stabilization to the same degree greater than or equal to 4, they are related by labelled 1-isotopy, meaning 1-isotopy that preserves the labelling consistently with the Wirtinger relations, and by the covering moves $c_1$ and $c_2$ in Figure* 24.

*Proof.* The first part of the statement is [**16**, Theorem 6] (see [**3**, Section 3] for a different proof based on Kirby diagrams), while the second part is [**3**, Theorem 1]. □

We remark that the orientability of a four-dimensional 2-handlebody does not imply the orientability of the labelled ribbon surfaces representing it as a simple branched covering of $B^4$. Nevertheless, by using the covering moves $c_1$ and $c_2$, any such labelled ribbon surface can be transformed into an orientable one, representing the same handlebody up to 2-deformations. In fact, those moves together with stabilization will enable us to easily convert any labelled planar diagram into a flat one.

As we anticipated at the end of Section 2, the covering moves $c_1$ and $c_2$ will also play a crucial role in the interpretation of Proposition 6 in terms of labelled flat planar diagrams. In this context, we can still use the moves $s_5$ and $s_6$ in Figure 6 and the flat isotopy moves of Figure 7, but not the isotopy moves of Figure 8 that involve half-twists.

On the other hand, the 1-isotopy moves of Figure 3, as well as the covering moves $c_1$ and $c_2$ themselves, which arise as three-dimensional moves, can also be thought of as moves of flat planar diagrams due to their flat presentation. However, when doing so one has to be careful to use them only accordingly to such fixed flat presentation.

Finally, (de)stabilization makes sense also for flat planar diagrams, being realizable in the elementary case as addition/deletion of a separate flat disc labelled $(i\ d+1)$ with $i \leqslant d$, to a $\Sigma_d$-labelled flat diagram. We call such a modification a *(de)stabilization move*.

As we will see shortly, in the presence of covering moves and stabilization all the moves of Figure 8 can be replaced by a unique move of flat planar diagrams (cf. Figure 27). But first we need the following lemma (cf. [**3**, **19**]).

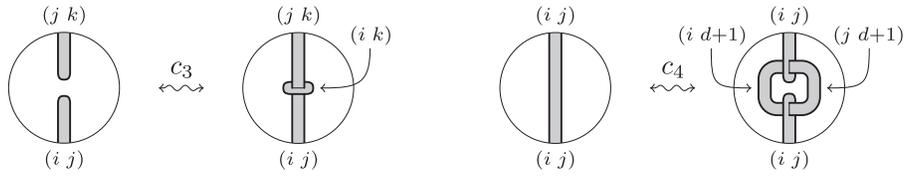

FIGURE 25. *Joining and splitting bands.*

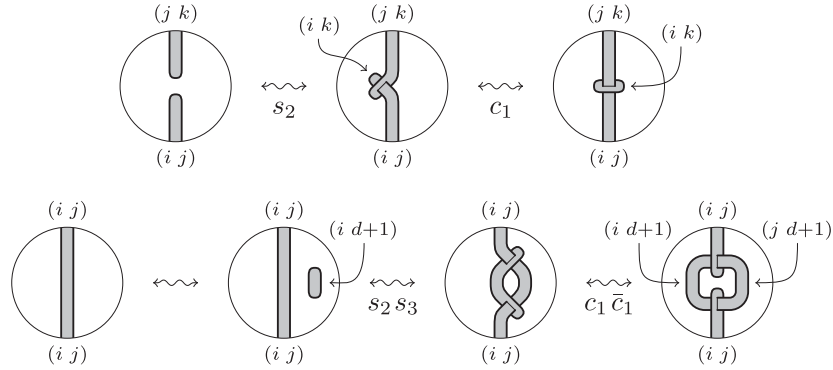

FIGURE 26. *Generating $c_3$ and $c_4$.*

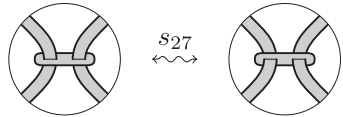

FIGURE 27. *The reversing move.*

LEMMA 7. *The covering moves $c_1$ and $c_2$ generate their symmetric $\bar c_1$ and $\bar c_2$, where terminal ribbon intersections of type* (h) *replace those of type* (g) *occurring in $c_1$ and $c_2$ (cf. Figure 4), modulo flat isotopy moves in Figure 7. Moreover, $c_1$ and $\bar c_1$ generate the self-symmetric moves $c_3$ and $c_4$ in Figure 25, where the left side of $c_4$ is assumed to be labelled in $\Sigma_d$, modulo the flat isotopy moves in Figure 7, the 1-isotopy moves $s_2$ and $s_3$ in Figure 3 and stabilization (actually required only for $c_4$).*

*Proof.* Move $\bar c_1$ can be reduced to the original move $c_1$, by applying the flat isotopy moves $\bar s_{18}$ and $\bar s_{17}$, respectively, to the left side and to right side. Similarly, move $\bar c_2$ can be reduced to $c_2$ modulo the flat isotopy moves $\bar s_{18}$ and $s_{11}$. Figure 26 shows how to generate moves $c_3$ and $c_4$. For the convenience of the reader, here and in the following figures, we indicate under the arrows the corresponding moves, omitting the flat isotopy ones. □

PROPOSITION 8. *Up to 2-deformations, any connected four-dimensional 2-handlebody can be represented as a simple branched covering of $B^4$ by a labelled flat planar diagram (cf. [3, Remark 2.7]). Two labelled flat planar diagrams represent 2-equivalent four-dimensional 2-handlebodies if and only if they are related by the* (de)*stabilization moves, the moves $s_1$ to $s_{18}$ and $s_{27}$ in Figures 3, 6, 7 and 27, and the moves $c_1$ and $c_2$ in Figure 24.*

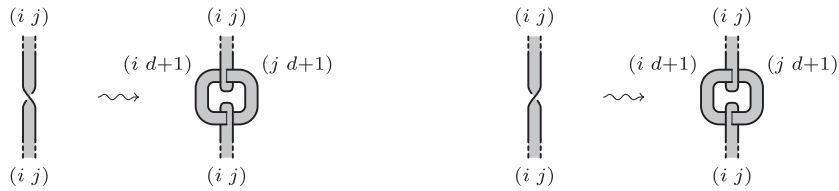

FIGURE 28. *Replacing half-twists.*

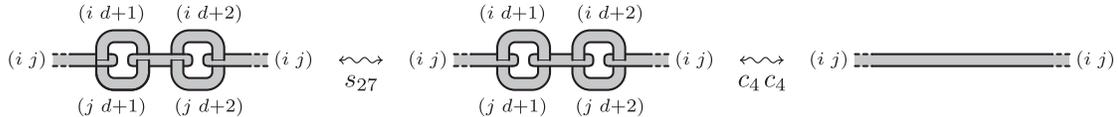

FIGURE 29. *Deriving move $s_{19}$.*

*Proof.* Proposition 6 tells us that any connected four-dimensional 2-handlebody can be represented by a labelled planar diagram. This can be made flat by replacing one by one in turn all the half-twist occurring in it as indicated in Figure 28, where $d$ is the degree of the covering. Of course, the degree of the covering increases by one at a single replacement; hence, we have different degrees $d$ when replacing different half-twists. Then, the final degree depends on the number of the half-twists (a different flattening procedure, which does not increase the degree, is described in [**3**, Remark 2.7]), while the final labelling depends on the order of the replacements and on the choice of $i$ (instead of $j$) for each one of them. In any case, we obtain a labelled flat planar diagram representing the same four-dimensional 2-handlebody as the original diagram, since each replacement can be thought as a move $c_4$ followed by a move $s_{23}$ or $\bar{s}_{23}$. This proves the first part of the proposition.

Now, assume that we have two labelled flat planar diagrams representing the same four-dimensional 2-handlebody up to 2-deformations. Then, by Proposition 6 they are related by a sequence of (de)stabilization moves, 1-isotopy moves $s_1$ to $s_{26}$ and covering moves $c_1$ and $c_2$. At each step of the sequence, if some half-twist is created by one of the moves $s_{19}$ to $s_{26}$, we replace it as described above. Then, we let the replacing configuration follow the original half-twist under the subsequent moves, until it disappears by the effect of one of the moves $s_{19}$ to $s_{26}$ again. In this way, we get a sequence of flat planar diagrams between the given ones, each related to the previous by the same move as in the original sequence, except that instead of the moves $s_{19}$ to $s_{26}$ we have their flat versions deriving from the replacement of half-twists. Then, to prove the second part of the proposition, it suffices to derive those flat versions from the moves prescribed in the statement. In doing that, we can also use the moves $\bar{c}_1$, $\bar{c}_2$, $c_3$ and $c_4$, thanks to Lemma 7, and all the symmetric moves $\bar{s}_5$ to $\bar{s}_{18}$, according to the discussion preceding Proposition 1. Moreover, since the labelling resulting from different choices in replacing the half-twists can be easily seen to be equivalent up to some moves $c_4$ and $s_{27}$ (possibly after renumbering the sheets of the covering), we can always assume it to be the most convenient one.

Move $s_{19}$ is realized in Figure 29, with one move $s_{27}$ and two moves $c_4$. Moves $s_{21}$ and $s_{22}$ can be derived in a similar way, with the help of some flat isotopy moves. Move $s_{20}$ can be skipped, being a consequence of $s_7$, $s_{23}$ and $s_{24}$ (in the special case when the bottom band is terminal), as we have already noted.

Moves $s_{23}$ and $s_{25}$ are obtained in Figure 30, by using moves $c_3$ and $s_2$ for the former and moves $s_{27}$ and $c_4$ for the latter.

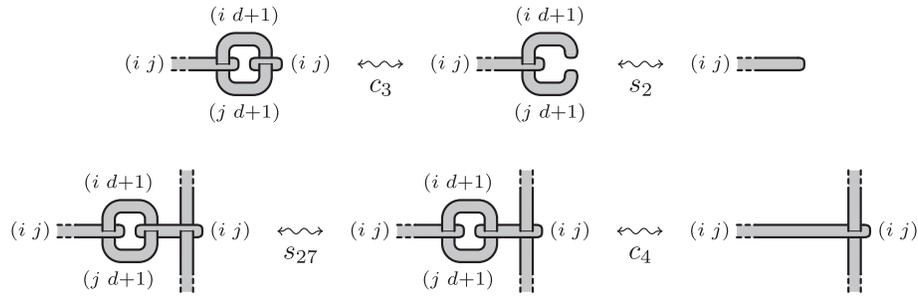

FIGURE 30. *Deriving moves $s_{23}$ and $s_{25}$.*

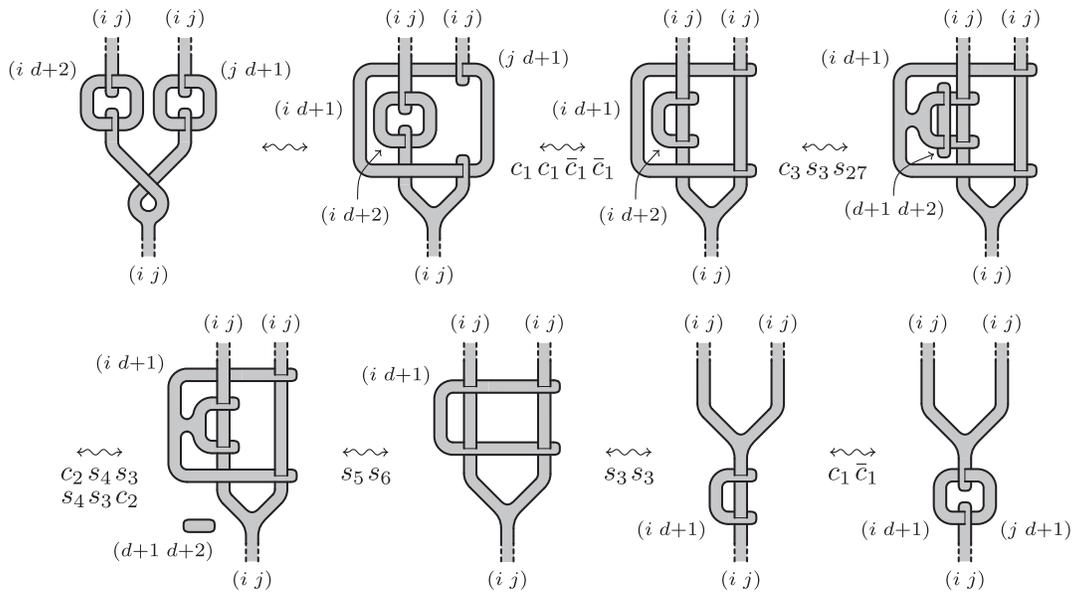

FIGURE 31. *Deriving move $s_{24}$.*

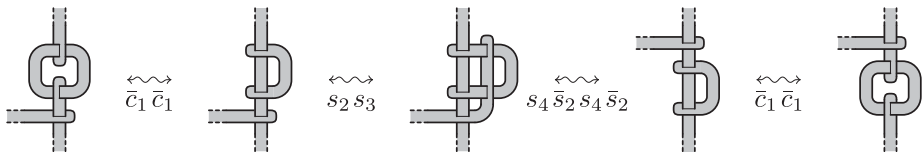

FIGURE 32. *Deriving moves $s_{26}$.*

Move $s_{24}$ and $s_{26}$ are considered, in an equivalent form up to $s_{19}$, in Figures 31 and 32, respectively. In particular, in Figure 32 we use the symmetric move $\bar{s}_2$, which can be easily reduced to $s_2$ modulo $s_5, s_{11}, s_{23}$ and $s_{25}$. It is worth remarking that move $s_{26}$ could also be realized by an obvious 1-isotopy, without involving the covering moves $c_1$ and $\bar{c}_1$, but this would require different planar projections of the 1-isotopy moves $s_3$ and $s_4$, much more difficult to get than $\bar{s}_2$. □

## 6. Lefschetz fibrations over $B^2$

A smooth map $f: W \to B^2$, with $W$ a smooth oriented compact four-manifold (possibly with corners), is called a *Lefschetz fibration* if the following properties hold.

(1) $f$ has a finite set $A = \{a_1, a_2, \ldots, a_n\} \subset \operatorname{Int} B^2$ of *singular values* and the restriction $f_|: W - f^{-1}(A) \to B^2 - A$ is a locally trivial fibre bundle, whose fibre is a compact connected orientable surface $F$ with (possibly empty) boundary, called the *regular fibre* of $f$.

(2) For any $a_i \in A$ the *singular fibre* $F_{a_i} = f^{-1}(a_i)$ contains only one singular point $w_i \in F_{a_i} \cap \operatorname{Int} W$ and there are local complex coordinates $(z_1, z_2)$ of $W$ and $z$ of $B^2$ centred at $w_i$ and $a_i$, respectively, such that $f: (z_1, z_2) \mapsto z = z_1^2 + z_2^2$.

If such coordinates $(z_1, z_2)$ can be chosen to preserve orientations (no matter whether $z$ does as well or not), then we call $w_i$ a *positive singular point*, otherwise we call it a *negative singular point*. Obviously, at a negative singular point we can always choose orientation preserving complex coordinates $(z_1, z_2)$ such that $f: (z_1, z_2) \mapsto z = z_1^2 + \bar{z}_2^2$.

Two Lefschetz fibrations $f: W \to B^2$ and $f': W' \to B^2$ are said to be *fibred equivalent* if there are orientation preserving diffeomorphisms $\varphi: B^2 \to B^2$ and $\widetilde{\varphi}: W \to W'$ such that $\varphi \circ f = f' \circ \widetilde{\varphi}$. Of course, in this case $\varphi$ restricts to a bijection $\varphi_|: A \to A'$ between the sets of singular values of $f$ and $f'$, respectively, while $\widetilde{\varphi}$ sends each singular point $w_i$ of $f$ into a singular point $w'_j$ of $f'$ with the same sign.

Note that the locally trivial fibre bundle $f_|$ in the definition of Lefschetz fibration $f: W \to B^2$ is oriented. Indeed, each regular fibre $F_x = f^{-1}(x) \cong F$ with $x \in B^2 - A$ has a preferred orientation, determined by the following rule: the orientation of $W$ at any point of $F_x$ coincides with the product of the orientation induced by the standard one of $B^2$ on any smooth local section of $f$ with the preferred one of $F_x$ in that order. In what follows, we will consider $F = F_* = f^{-1}(*)$ endowed with this preferred orientation, for the fixed base point $* \in S^1$.

On the other hand, any singular fibre $F_{a_i}$ is an orientable surface away from the singular point $w_i$ and the preferred orientation of the regular fibres coherently extends to $F_{a_i} - \{w_i\}$. Moreover, when $\operatorname{Bd} F \neq \emptyset$, by putting $\operatorname{Bd} F_{a_i} = \operatorname{Bd}(F_{a_i} - \{w_i\})$ for $i = 1, 2, \ldots, n$ and $T = \bigcup_{x \in B^2} \operatorname{Bd} F_x \subset \operatorname{Bd} W$, we have that $f_{|T}: T \to B^2$ is a trivial bundle with fibre $\operatorname{Bd} F$. In this case, corners naturally occur along $T \cap f^{-1}(S^1) = \bigcup_{x \in S^1} \operatorname{Bd} F_x \subset \operatorname{Bd} W$.

The structure of $f$ over a small disc $D_i$ centred at a singular value $a_i$ is given by the following commutative diagram, where: $\gamma^{\pm 1}: F \to F$ is a Dehn twist along a cycle $c \subset F$, and it is positive (right-handed) or negative (left-handed) according to the sign of the singular point $w_i \in F_{a_i}$; $T(\gamma^\pm) = F \times [0,1]/((\gamma^{\pm 1}(x), 0) \sim (x, 1) \,\forall x \in F)$ is the mapping torus of $\gamma^{\pm 1}$ and $\pi: T(\gamma^{\pm 1}) \to S^1 \cong [0,1]/(0 \sim 1)$ is the canonical projection; the singular fibre $F_{a_i} \cong F/c$ has a node singularity at $w_i$, which is positive or negative according to the sign of $w_i$; $\varphi$ and $\widetilde{\varphi}$ are orientation preserving diffeomorphisms such that the cycles $c_x = \widetilde{\varphi}([c,s], t) \subset F_x$ collapse to $w_i$ as $x = \varphi(s, t) \to a_i$. (cf. [**9**] or [**13**])

$$\begin{array}{ccccccc}
T(\gamma^\pm) \times (0,1] & \xrightarrow{\widetilde{\varphi}} & f^{-1}(D_i) - F_{a_i} & \subset & f^{-1}(D_i) & \supset & F_{a_i} \\
\downarrow {\scriptstyle \pi \times \operatorname{id}} & & \downarrow & & \downarrow {\scriptstyle f_|} & & \downarrow \\
S^1 \times (0,1] & \xrightarrow{\varphi} & D_i - \{a_i\} & \subset & D_i & \ni & a_i
\end{array}$$

Because of this collapsing, the cycles $c_x$ are called *vanishing cycles*. We point out that they are well defined up to ambient isotopy of the fibres $F_x$, while the cycle $c \subset F$ is only defined up to diffeomorphisms of $F$, depending on the specific identification $F \cong F_x$ induced by $\widetilde{\varphi}$. The indeterminacy of the cycle $c \subset F$ can be resolved if a Hurwitz system $(\alpha_1, \alpha_2, \ldots, \alpha_n)$ for $A$ is given. In fact, we can choose $\widetilde{\varphi}$ such that the induced identification $F \cong F_x$ coincides with the one deriving from any trivialization of $f_|$ over any arc $\alpha'_i$ joining $*$ to $x$ in $B^2 - A$ and running along $\alpha_i$ outside $D_i$. In this way, we get a cycle $c_i \subset F$ well defined up to ambient isotopy of $F$, which represents the vanishing cycles at $w_i$ in the regular fibre $F = F_*$. We denote by $\gamma_i$, the positive (right-handed) or negative (left-handed) Dehn twists of $F$ along $c_i$ corresponding to $\gamma^{\pm 1}$ in the above diagram, which is uniquely determined up to ambient isotopy of $F$ as well. We call $c_i$ the *vanishing cycle* of $f$ over $a_i$ and $\gamma_i$ the *mapping monodromy* of $f$ over $a_i$ (with respect to the given Hurwitz system).

The Lefschetz fibration $f: W \to B^2$ with the set of singular values $A \subset B^2$ turns out to be uniquely determined, up to fibred equivalence, by its mapping monodromy sequence $(\gamma_1, \gamma_2, \ldots, \gamma_n)$ with respect to any given Hurwitz system $(\alpha_1, \alpha_2, \ldots, \alpha_n)$ for $A$. Of course, we can identify $F$ with the standard compact connected oriented surface $F_{g,b}$ with genus $g \geqslant 0$ and $b \geqslant 0$ boundary components, and think of each $\gamma_i$ as a Dehn twist of $F_{g,b}$. Actually, we will represent them as signed cycles in $F_{g,b}$.

According to our discussion about Hurwitz systems in Section 3, mapping monodromy sequences associated to different Hurwitz systems are related by elementary transformations, changing a given sequence of Dehn twists $(\gamma_1, \gamma_2, \ldots, \gamma_n)$ into $(\gamma'_1, \gamma'_2, \ldots, \gamma'_n)$ with $\gamma'_i = \gamma_i \gamma_{i+1} \gamma_i^{-1}$, $\gamma'_{i+1} = \gamma_i$ and $\gamma'_k = \gamma_k$ for $k \neq i, i+1$, for some $i < n$, and their inverses. We call this transformation the *twist sliding* of $\gamma_{i+1}$ over $\gamma_i$, and its inverse the twist sliding of $\gamma'_i$ over $\gamma'_{i+1}$.

When considering $f$ up to fibred equivalence, we can always assume $(\alpha_1, \alpha_2, \ldots, \alpha_n)$ to be the standard Hurwitz system and consider $(\gamma_1, \gamma_2, \ldots, \gamma_n)$ as an abstract sequence of Dehn twists of $F_{g,b}$ without any reference to a specific Hurwitz system. In this perspective, twist slidings can be interpreted as fibred isotopy moves, and two sequences of Dehn twists of $F_{g,b}$ represent fibred equivalent Lefschetz fibrations if and only if they are related by: (1) the simultaneous action of $\mathcal{M}_{g,b} = \mathcal{M}_+(F_{g,b})$ on the vanishing cycles, where $\mathcal{M}_+$ denotes the positive mapping class group consisting of all isotopy classes of orientation preserving diffeomorphisms fixing the boundary, to take into account possibly different identifications $F \cong F_{g,b}$; (2) twist slidings (hence cyclic shift of the $\gamma_i$'s as well), to pass from one Hurwitz system to another.

Actually, any sequence $(\gamma_1, \gamma_2, \ldots, \gamma_n)$ of positive or negative Dehn twists of $F_{g,b}$ does represent in this way a Lefschetz fibration $f: W \to B^2$ with regular fibre $F \cong F_{g,b}$, uniquely determined up to fibred equivalence. Such a Lefschetz fibration $f$ can be constructed as described below.

The most elementary non-trivial Lefschetz fibrations over $B^2$ are the *Hopf fibrations* $h_\pm : H_\pm \to B^2$. These are defined as $h_+(z_1, z_2) = z_1^2 + z_2^2$ and $h_-(z_1, z_2) = z_1^2 + \bar{z}_2^2$ for all $(z_1, z_2) \in \mathbb{C}^2$, and $H_\pm = \{(z_1, z_2) \in \mathbb{C} \mid |z_1|^2 + |z_2|^2 \leqslant 2 \text{ and } |h_\pm(z_1, z_2)| \leqslant 1\} \cong B^4$ (up to smoothing the corners). Their regular fibre is an annulus $F \cong F_{0,2}$ and they have $w_1 = (0,0)$ and $a_1 = 0$ as the unique singular point and singular value respectively, while the mapping monodromy $\gamma_1$ is the right-handed Dehn twist along the unique vanishing cycle $c_1$ represented by the core of $F$ for $h_+$, and the left-handed Dehn twist along $c_1$ for $h_-$. Furthermore, for each $x \in S^1$, the regular fibre $F_x = h_+^{-1}(x)$ or $h_-^{-1}(x)$ forms a left- or right-handed, respectively, full twist as an embedded closed band in $\operatorname{Bd} H_\pm \cong S^3$. We call $h_+$ and $h_-$, the *positive* and *negative*, respectively, Hopf fibration.

Now we introduce a fibre gluing operation, which will allow us to build up any other non-trivial Lefschetz fibration over $B^2$ by using Hopf fibrations as the basic blocks, and to describe the equivalence moves in Section 7 as well.

Let $f_1: W_1 \to B^2$ and $f_2: W_2 \to B^2$ be two Lefschetz fibrations with regular fibres $F_1 = f_1^{-1}(*)$ and $F_2 = f_2^{-1}(*)$, respectively, and let $\eta: G_1 \to G_2$ be a diffeomorphism between two smooth subsurfaces (possibly with corners) $G_1 \subset F_1$ and $G_2 \subset F_2$ such that $F = F_1 \cup_\eta F_2 = (F_1 \sqcup F_2)/(x \sim \eta(x) \, \forall x \in G_1)$ is a smooth surface (possibly with corners). For $i = 1, 2$, we can consider $F_i \subset F$ and hence $\operatorname{Fr}_F F_i \subset \operatorname{Bd} F_i$. Moreover, once a trivialization $\varphi_i: T_i = \bigcup_{x \in B^2} \operatorname{Bd} f_i^{-1}(x) \to B^2 \times \operatorname{Bd} F_i$ of the restriction $f_{i|}: T_i \to B^2$ is chosen such that $\varphi_i(x) = (*, x)$ for all $x \in \operatorname{Bd} F_i$, we can extend $f_i$ to a Lefschetz fibration $\widehat{f}_i: \widehat{W}_i \to B^2$ with $F$ as the regular fibre, in the following way. We put $\widehat{W}_i = W_i \cup_{\varphi'_i} (B^2 \times \operatorname{Cl}_F(F - F_i))$, where $\varphi'_i$ is the restriction of $\varphi_i$ to $T'_i = \varphi_i^{-1}(B^2 \times \operatorname{Fr}_F F_i) \subset T_i$, and define $\widehat{f}_i$ to coincide with the projection onto the first factor in $B^2 \times \operatorname{Cl}_F(F - F_i)$. Then, let $I_1, I_2 \subset S^1$ denote two intervals, respectively, ending to and starting from $* \in S^1$ (in the counterclockwise orientation) and let $\psi_i \colon \widehat{f}_i^{-1}(I_i) \to I_i \times F$ be any trivialization of the restrictions $\widehat{f}_i| \colon \widehat{f}_i^{-1}(I_i) \to I_i$ such that $\psi_i(x) = (*, x)$ for all $x \in F$, with $i = 1, 2$. Finally, we define a new Lefschetz fibration $f_1 \#_\eta f_2 \colon W_1 \#_\eta W_2 \to B^2 \# B^2 \cong B^2$, where $B^2 \# B^2 = B^2 \cup_\rho B^2$ is the boundary connected sum given by an orientation reversing identification $\rho \colon I_1 \to I_2$, by putting $W_1 \#_\eta W_2 = \widehat{W}_1 \cup_\psi \widehat{W}_2$ and $f_1 \#_\eta f_2 = \widehat{f}_1 \cup_\psi \widehat{f}_2$, with $\psi = \psi_2^{-1} \circ (\rho \times \operatorname{id}_F) \circ \psi_1 \colon f_1^{-1}(I_1) \to f_2^{-1}(I_2)$. A straightforward verification shows that, this is well defined up to fibred equivalence, depending only on $f_1$, $f_2$ and $\eta$, but not on the various choices involved in its construction.

We call the Lefschetz fibration $f_1 \#_\eta f_2 \colon W_1 \#_\eta W_2 \to B^2$ the *fibre gluing* of $f_1$ and $f_2$ through the diffeomorphism $\eta \colon G_1 \to G_2$. It has regular fibre $F = F_1 \cup_\eta F_2$. Moreover, under the identification $B^2 \# B^2 \cong B^2$, its set of singular values is the disjoint union $A = A_1 \sqcup A_2 \subset B^2$ of those of $f_1$ and $f_2$, while a mapping monodromy sequence for it is given by the juxtaposition of two given sequences for $f_1$ and $f_2$, with all the Dehn twists thought of as acting on $F$, through the inclusions $F_i \subset F$.

Up to fibred equivalence, fibre gluing is weakly associative in the sense that the equivalence between $(f_1 \#_{\eta_1} f_2) \#_{\eta_2} f_3$ and $f_1 \#_{\eta_1} (f_2 \#_{\eta_2} f_3)$ holds under the assumption (not always true) that all the gluings appearing in both the expressions make sense. This fact easily follows from the definition and allows us to write $f_1 \#_{\eta_1} f_2 \#_{\eta_2} \cdots \#_{\eta_{n-1}} f_n$ without brackets. Still up to fibred equivalence, fibre gluing is also commutative, being the monodromy sequences for $f_1 \#_\eta f_2$ and $f_2 \#_{\eta^{-1}} f_1$ related by a cyclic shift.

Finally, it is worth noticing that the fibre gluing $f_1 \#_\eta f_2$ reduces to the usual fibre sum (cf. [**9**]) when $G_i = F_i = F$ for $i = 1, 2$ and $\eta = \operatorname{id}_F$. On the other hand, as we will see in the next section, it also includes as a special case the Hopf plumbing.

Then, given any sequence of positive or negative Dehn twists $(\gamma_1, \gamma_2, \ldots, \gamma_n)$ of $F_{g,b}$, a Lefschetz fibration $f \colon W \to B^2$ with that mapping monodromy sequence is provided by the fibre gluing $f = f_0 \#_{\eta_1} h_1 \#_{\eta_2} h_2 \cdots \#_{\eta_n} h_n$, where: $f_0$ is the product fibration $B^2 \times F_{g,b} \to B^2$; $h_i$ is a positive or negative Hopf fibration according to the sign of $\gamma_i$; $\eta_i \colon N_i \to F_i$ is an orientation preserving diffeomorphism between a regular neighbourhood $N_i \subset \operatorname{Int} F_{g,b}$ of the cycle $c_i$ along which $\gamma_i$ occurs and the regular fibre $F_i$ of $h_i$.

The total space $W$ of any Lefschetz fibration $f \colon W \to B^2$ has a natural four-dimensional 2-handlebody structure $H_f$, induced by $\|f\|^2 \colon W \to [0, 1]$ as a Morse function away from 0 (see [**9**] or [**13**]). For our aims, it is more convenient to derive such handlebody structure of $W$ from a mapping monodromy sequence representing $f$, through the corresponding fibre gluing decomposition. This is the point of view adopted in the next proposition.

PROPOSITION 9. *Any Lefschetz fibration $f \colon W \to B^2$ determines a four-dimensional 2-handlebody decomposition $H_f$ of $W$, well defined up to 2-deformations. Moreover, the 2-equivalence class of $H_f$ is invariant under fibred equivalence of Lefschetz fibrations.*

*Proof.* Let $(\gamma_1, \gamma_2, \ldots, \gamma_n)$ be the mapping monodromy sequence of $f$ associated to any Hurwitz system for the set of singular values $A = \{a_1, a_2, \ldots, a_n\} \subset B^2$. Then, a four-dimensional 2-handlebody decomposition of $W$ based on the corresponding fibre gluing presentation $f = f_0 \#_{\eta_1} h_1 \#_{\eta_2} h_2 \cdots \#_{\eta_n} h_n$ described above, can be constructed as follows.

We start with any handlebody decomposition $H_F$ of $F_{g,b}$ and consider the induced four-dimensional 2-handlebody decomposition $B^2 \times H_F$ of the product $B^2 \times F_{g,b}$ (actually, this can be assumed to have no 2-handles if $b > 0$). Then, for all $i = 1, 2, \ldots, n$, we define the 2-handle $H_i^2$ as the total space $H_\pm \cong B^4$ of the Hopf fibration $h_i \colon H_\pm \to B^2$, attached to $B^2 \times F_{g,b}$ through the map $(\rho \times \eta_i)^{-1}: I_2 \times F_i \to I_1 \times F_{g,b}$. The attaching sphere of $H_i^2$ is a copy $\{y_i\} \times c_i \subset S^1 \times F_{g,b}$ of the vanishing cycle $c_i \subset F_{g,b}$, while its attaching framing turns out to be $\pm 1$ with respect to the one given by $F_{g,b}$. Namely, we have $-1$ or $+1$ if the singular point $w_i$, hence the Hopf fibration $h_i$, is positive or negative, respectively, due to the full twist with the opposite sign formed by the fibre of $h_i$ in $\text{Bd}\, H_\pm \cong S^3$. The points $y_1, y_2, \ldots, y_n$ are ordered along $S^1 - \{*\}$ according to the counterclockwise orientation.

Now, we let $H_f = (B^2 \times H_F) \cup H_1^2 \cup H_2^2 \cup \cdots \cup H_n^2$ be the handlebody decomposition of $W$ just constructed and observe that, up to handle isotopy, it depends only on the choices of the mapping monodromy sequence $(\gamma_1, \gamma_2, \ldots, \gamma_n)$ and of the handlebody decomposition $H_F$ of $F_{g,b}$.

The well definedness of $H_f$ up to 2-deformations and its invariance under fibred isotopy of $f$ are immediate consequences of the following facts: (1) diffeomorphic two-dimensional handlebodies are always 2-equivalent; (2) any twist sliding in the mapping monodromy sequence $(\gamma_1, \gamma_2, \ldots, \gamma_n)$ induces a number of 2-handle slidings on the handlebody $H_f$, one for each transversal intersection between the cycles involved in the twist sliding. $\square$

As proved by Harer in his thesis [**11**], up to 2-equivalence any four-dimensional 2-handlebody decomposition of $W$ can be represented by a Lefschetz fibration $f: W \to B^2$ according to Proposition 9. This could also be derived from Proposition 8, by applying the braiding procedure discussed in Section 4 to the labelled flat diagram representing the given handlebody decomposition (cf. Proposition 10).

The natural question of how to relate any two such representations of 2-equivalent four-dimensional 2-handlebodies will be answered in the next sections, by using the branched covering representation of Lefschetz fibrations we are going to describe in the final part of this section. As a preliminary step, let us briefly discuss the notion of mapping monodromy homomorphism of a Lefschetz fibration.

Let $f: W \to B^2$ be a Lefschetz fibration with set of singular values $A = \{a_1, a_2, \ldots, a_n\} \subset B^2$ and regular fibre $F \cong F_{g,b}$. Since $\pi_1(B^2 - A)$ is freely generated by any Hurwitz system $(\alpha_1, \alpha_2, \ldots, \alpha_n)$ for $A$, the corresponding mapping monodromy sequence $(\gamma_1, \gamma_2, \ldots, \gamma_n)$ for $f$ gives rise to a homomorphism $\omega_f: \pi_1(B^2 - A) \to \mathcal{M}_{g,b} = \mathcal{M}_+(F_{g,b}) \cong \mathcal{M}_+(F)$ such that $\omega_f(\alpha_i) = \gamma_i$ for every $i = 1, 2, \ldots, n$ (remember that we denote by $\mathcal{M}_+$ the positive mapping class group, consisting of all isotopy classes of orientation preserving diffeomorphisms fixing the boundary).

We call $\omega_f: \pi_1(B^2 - A) \to \mathcal{M}_{g,b}$ the *mapping monodromy* of $f$. Note that $\omega_f$ is defined only up to conjugation in $\mathcal{M}_{g,b}$, depending on the chosen identification $F \cong F_{g,b}$. On the contrary, $\omega_f$ does not depend on the specific sequence $(\gamma_1, \gamma_2, \ldots, \gamma_n)$, admitting an intrinsic definition not based on the choice of a Hurwitz system.

We outline this definition of $\omega_f$, to emphasize that for $\text{Bd}\, F \neq \emptyset$ (that is $b > 0$) it also involves the choice of a trivialization $\varphi: T \cong B^2 \times \text{Bd}\, F$ of the bundle $f_{|T}: T \to B^2$, such that $\varphi(x) = (*, x)$ for all $x \in \text{Bd}\, F$. This is used to achieve the condition that $\omega_f([\lambda])$ fixes $\text{Bd}\, F$ for any $[\lambda] \in \pi_1(B^2 - A)$, as follows. Given the loop $\lambda: [0, 1] \to B^2 - A$, we first consider the commutative diagram below, where the total space of the induced fibre bundle $\lambda^*(f_|)$ is identified with $[0, 1] \times F$ by a trivialization of $\lambda^*(f_|)$, in such a way that $\widetilde{\lambda}(0, x) = x$ for all $x \in F$ and $\widetilde{\lambda}(t, x) = \varphi(\lambda(t), x)$ for all $x \in \text{Bd}\, F$.

$$
\begin{array}{ccc}
[0,1] \times F & \xrightarrow{\widetilde{\lambda}} & W - f^{-1}(A) \\
{\scriptstyle \lambda^*(f_|)} \downarrow & & \downarrow {\scriptstyle f_|} \\
[0,1] & \xrightarrow{\lambda} & B^2 - A
\end{array}
$$

Then we put $\omega_f([\lambda]) = [\widetilde{\lambda}_1]$, with $\widetilde{\lambda}_1 : F \to F$ the diffeomorphism fixing $\operatorname{Bd} F$, defined by $\widetilde{\lambda}_1(x) = \widetilde{\lambda}(1, x)$ for all $x \in F$. It is not difficult to see that this definition is independent on the specific choice of $\varphi$, since different choices are fibrewise isotopic.

From the classical theory of fibre bundles, we know that $\omega_f$ uniquely determines the restriction of $f_|$ over $B^2 - A$ up to fibred equivalence. But, in general, it does not determine the whole Lefschetz fibration $f$. In fact, when considering $\gamma_i$ as an element of $\mathcal{M}_{g,b}$ the sign of it as a Dehn twist gets lost if the cycle $c_i$ is homotopically trivial in $F_{g,b}$, being in this case $\gamma_i$ and $\gamma_i^{-1}$ both isotopic to the identity, while there is no loss of information in the non-trivial cases. Therefore, the mapping monodromy $\omega_f : \pi_1(B^2 - A) \to \mathcal{M}_{g,b}$ uniquely determines the Lefschetz fibration $f : W \to B^2$ up to fibred equivalence only under the assumption that no vanishing cycle of $f$ is homotopically trivial in $F_{g,b}$, in which case $f$ is called a *relatively minimal* Lefschetz fibration. While in the presence of trivial vanishing cycles, $f$ turns out to be determined only up to blow-ups (cf. [**9**]).

A Lefschetz fibration $f : W \to B^2$ is called *allowable* if its regular fibre $F$ has boundary $\operatorname{Bd} F \neq \emptyset$ and all its vanishing cycles are homologically non-trivial in $F$ (hence $f$ is relatively minimal as well). Of course, since the property of being homologically non-trivial is invariant under the action of $\mathcal{M}_+(F)$, it is enough to verify it for the vanishing cycles $c_i$ of any given monodromy sequence $(\gamma_1, \gamma_2, \ldots, \gamma_n)$ for $f$.

According to Loi and Piergallini [**15**] and Zuddas [**23**], allowable Lefschetz fibrations are closely related to simple covering of $B^2 \times B^2$ branched over braided surfaces. This relation is established by the next two propositions. But first we recall the notion of lifting braids (cf. [**2**, **17**] or [**23**]).

Let $q : F \to B^2$ be a simple covering branched over the finite set $\widetilde{*} \in \Gamma_m B^2$. Then a braid $\beta \in \mathcal{B}_m \cong \pi_1(\Gamma_m B^2, \widetilde{*})$ is said to be *liftable* with respect to $q$, when such is the terminal diffeomorphism $h_1 : (B^2, \widetilde{*}) \to (B^2, \widetilde{*})$ of an ambient isotopy $(h_t)_{t \in [0,1]}$ of $B^2$ that fixes $S^1$ and realizes $\beta$ as the loop $t \mapsto h_t(\widetilde{*})$ in $\Gamma_m B^2$, meaning that there exists a diffeomorphism $\widetilde{h}_1 : F \to F$ such that $q \circ \widetilde{h}_1 = h_1 \circ q$. We denote by $\mathcal{L}_q \subset \mathcal{B}_m$ the subgroup of the liftable braids and by $\lambda_q : \mathcal{L}_q \to \mathcal{M}_+(F)$ the *lifting homomorphism*, that sends $\beta$ to the isotopy class of $\widetilde{h}_1$. It turns out that, if $\beta \in \mathcal{L}_q$ is a positive or negative half-twist around an arc $b \subset B^2$, then $\lambda_q(\beta)$ is, respectively, the positive or negative Dehn twist along the unique cycle component of $q^{-1}(b) \subset F$. Actually, every compact connected orientable surface $F$ with boundary $\operatorname{Bd} F \neq \emptyset$ admits a branched covering $q : F \to B^2$, such that any Dehn twist in $\mathcal{M}_+(F)$ along a homologically non-trivial cycle of $F$ can be represented in this way. This result dates back to the 1970s in the special case when $\operatorname{Bd} F$ is connected, while it was proved in [**23**] in the general case.

PROPOSITION 10. *Let $p : W \to B^2 \times B^2$ be a simple covering branched over a braided surface $S \subset B^2 \times B^2$. Then, the composition $f = \pi \circ p : W \to B^2$, where $\pi : B^2 \times B^2 \to B^2$ is the projection onto the first factor, is an allowable Lefschetz fibration. The set $A \subset B^2$ of singular values of $f$ coincides with the branch set of the branched covering $\pi_{|S} : S \to B^2$, and the mapping monodromy of $f$ is the lifting $\omega_f = \lambda_q \circ \omega_S$ of the braid monodromy $\omega_S$ of $S$, through the branched covering $q = p_| : F \cong f^{-1}(*) \to \pi^{-1}(*) \cong B^2$ representing the regular fibre $F$ of $f$ (note that $\operatorname{Im} \omega_S \subset \mathcal{L}_q$). Moreover, the four-dimensional 2-handlebody decompositions $H_f$ and $H_p$ of $W$, given by Propositions 9 and 5, respectively, coincide (up to 2-equivalence).*

*Proof.* The first part of the statement is a special case of Proposition 1 of [**15**]. The proof of the equations $\omega_f = \lambda_q \circ \omega_S$ and $H_f = H_p$ is just a matter of comparing the definitions. $\square$

PROPOSITION 11. *Any allowable Lefschetz fibration $f : W \to B^2$ factorizes as a composition $f = \pi \circ p$, where $\pi : B^2 \times B^2 \to B^2$ is the projection onto the first factor and $p : W \to B^2 \times B^2$* is a simple covering branched over a braided surface $S \subset B^2 \times B^2$ (actually, $p$ could be assumed to have degree 3 when $\mathrm{Bd}\, F$ is connected, but we will not need this fact here).

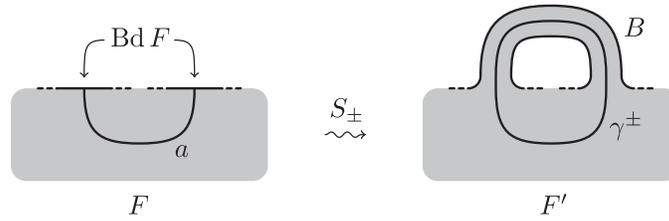

Figure 33. *Hopf stabilization.*

*Proof.* This is a special case of [**15**, Proposition 2]. □

In light of Propositions 10 and 11, up to composition with $\pi$, labelled braided surfaces $S \subset B^2 \times B^2$ (in fact, their band presentations) representing simple branched coverings $p\colon W \to B^2 \times B^2$, can be used to represent allowable Lefschetz fibrations $f\colon W \to B^2$ as well. Under this representation, labelled braided isotopy and band sliding for labelled braided surfaces, respectively, correspond to fibred equivalence and twist sliding for Lefschetz fibrations.

## 7. The 2-equivalence moves

In this section, we describe some operations on a Lefschetz fibration $f\colon W \to B^2$, which preserve the 2-equivalence class of the four-dimensional 2-handlebody decomposition $H_f$ induced on the total space $W$, hence the smooth topological type of $W$ as well.

*S move.* This is the well known Hopf stabilization (or plumbing) of a Lefschetz fibration with bounded regular fibre. In terms of fibre gluing, it can be defined as follows.

Let $f\colon W \to B^2$ be a Lefschetz fibration with regular fibre $F$ such that $\mathrm{Bd}\, F \neq \emptyset$, $a \subset F$ be a proper smooth arc and $G \subset F$ be a regular neighbourhood of $a$. On the other hand, let $F_h \cong F_{0,2}$ be the regular fibre of the Hopf fibration $h_\pm\colon H_\pm \to B^2$ and $G_h \subset F_h$ be a regular neighbourhood of a transversal arc in the annulus $F_h$. Then, the *positive* or *negative Hopf stabilization* of $f$ is the fibre gluing $f' = f \#_\eta h_+$ or $f' = f \#_\eta h_-$, respectively, with $\eta\colon G \to G_h$ a diffeomorphism such that $\eta(a)$ is an arc in the vanishing cycle $c_h$ of $h_\pm$ (the core of $F_h$).

Up to fibred equivalence, the stabilization $f \#_\eta h_\pm$ turns out to depend only on the fibred equivalence class of $f$ and on the isotopy class of $a$ in $F$. In fact, its regular fibre $F'$ is given by the attachment of a new band $B$ to $F$ along the arcs $G \cap \mathrm{Bd}\, F$, while a mapping monodromy sequence for $f \#_\eta h_\pm$ can be obtained from one for $f$, by inserting anywhere in the sequence a positive or negative Dehn twist $\gamma^\pm$ along a new vanishing cycle $c \subset F'$ running once over the band $B$ (see Figure 33).

In what follows, we will denote by $S_\pm\colon f \rightsquigarrow f' = f \#_\eta h_\pm$ the positive or negative *Hopf stabilization move* and by $S_\pm^{-1}\colon f' \rightsquigarrow f$ the *Hopf destabilization move* inverse of it. The latter can be performed on $f'$ whenever the regular fibre $F'$ has a 1-handle (in some handlebody decomposition of it) that is traversed once by only one vanishing cycle $c$. By a (de)stabilization of a Lefschetz fibration we mean the result of consecutive Hopf (de)stabilizations.

Note that, if the end points of the arc $a$ in the above definition belong to different components of $\mathrm{Bd}\, F$, then $F'$ has one less boundary component than $F$. Hence, the boundary of the regular fibre of a Lefschetz fibration can be made connected by a suitable sequence of Hopf stabilizations.

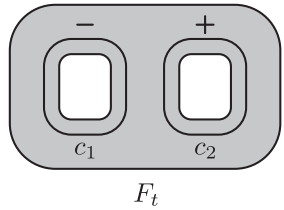

Figure 34. *The Lefschetz fibration t.*

Looking at the handlebody decomposition $H_f$, we see that a Hopf stabilization results into an addition of a cancelling pair of handles of indices 1 and 2. Namely, the 1-handle derives from the new band $B$, and a 2-handle is attached along a parallel copy of $c$ with framing $-1$ or $+1$ with respect to the fibre, depending on the stabilization being positive or negative, respectively. Thus, the 2-equivalence class of $H_f$ is preserved by Hopf stabilization.

For an allowable Lefschetz fibration $f: W \to B^2$ represented by a $\Sigma_d$-labelled braided surface $S$ according to Proposition 11, a Hopf (de)stabilization corresponds to an elementary labelled (de)stabilization of $S$ as defined in Section 3. This changes the fibration $f$, but not the covering $p: W \to B^2 \times B^2$ up to smooth equivalence (after smoothing the corners), and it should not be confused with the covering stabilization obtained by the addition of an extra separate sheet to $S$ with monodromy $(i\ d+1)$ for some $i \leqslant d$, which, on the contrary, changes the covering $p$, but not the Lefschetz fibration $f$ up to fibred equivalence.

As a consequence, we have that allowability of Lefschetz fibrations is preserved by Hopf (de)stabilization. Moreover, Proposition 3 implies that any allowable Lefschetz fibration admits a positive stabilization represented by a labelled braided surface with monotonic bands. Finally, it is worth mentioning that the $S$ move has been used by the third author in [24] to construct universal Lefschetz fibrations (the analogous of universal bundles).

*T move.* This is a new move, which corresponds to particular 2-deformations of the handlebody decomposition $H_f$ of the total space $W$ of a Lefschetz fibration $f: W \to B^2$. Like Hopf stabilization, the $T$ move applies only if the regular fibre $F$ of $f$ has non-empty boundary, but an extra condition is required on the mapping monodromy of $f$. This condition can be expressed by assuming that $f = f_0 \#_\eta t$, where $f_0: W_0 \to B^2$ is any Lefschetz fibration with bounded regular fibre, while the specific Lefschetz fibration $t: B^4 \to B^2$ and the gluing map $\eta$ are as follows.

Figure 34 describes $t: B^4 \to B^2$ in terms of its regular fibre $F_t \cong F_{0,3}$ and its monodromy sequence $(\gamma_1, \gamma_2)$. Here, the Dehn twists $\gamma_1$ and $\gamma_2$ are represented by the signed vanishing cycles $c_1$ and $c_2$, respectively, parallel to the inner boundary components. We assume the twists to have opposite signs, since this is enough for our purposes, but this assumption could be relaxed, as we will see later. Moreover, we note that $\gamma_1$ and $\gamma_2$ can be interchanged in the sequence, since $c_1$ and $c_2$ are disjoint.

The gluing map $\eta: G_0 \to G_t$ is shown in Figure 35. The surface $G_0 \subset F_0$ is an annulus in the regular fibre $F_0$ of $f_0$, whose core is the oriented cycle $a \subset \mathrm{Int}\, F_0$ and whose boundary meets $\mathrm{Bd}\, F_0$ along the four arcs indicated in the figure, in such a way that the oriented transversal arcs $r$ and $s$ are properly embedded in $F_0$. On the left side of the figure, we see the annulus $G_t \subset F_t$ with the oriented cycle and arcs corresponding to $a$, $r$ and $s$ under $\eta$. While the right side gives an analogous description of a different gluing map $\eta': G_0 \to G_t$. The outer boundary component of $G_t$ coincides with that of $F_t$, while the inner one meets along two arcs those of $F_t$. Of course, the given data uniquely determine the rest of $\eta$ and $\eta'$ up to isotopy.

We call a $T$ *move* the transformation $T: f \leftrightsquigarrow f'$, with $f = f_0 \#_\eta t: W = W_0 \#_\eta B^4$ and $f' = f_0 \#_{\eta'} t: W' = W_0 \#_{\eta'} B^4$. A more explicit description of such move is provided in Figure 36.

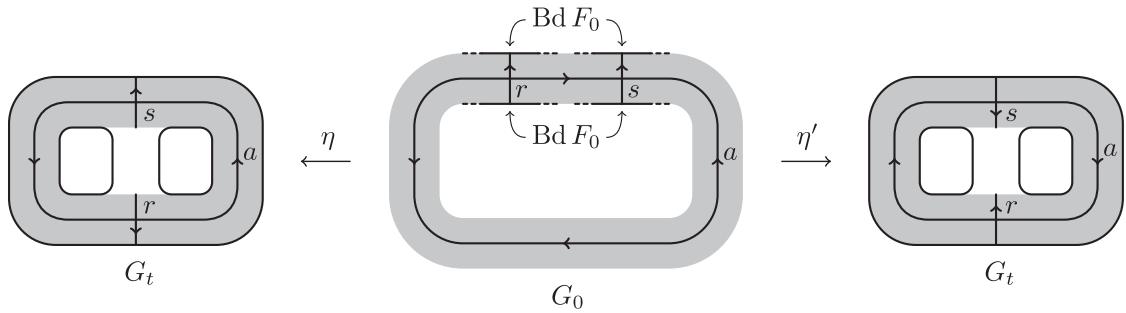

Figure 35. *The gluing maps between $f_0$ and $t$.*

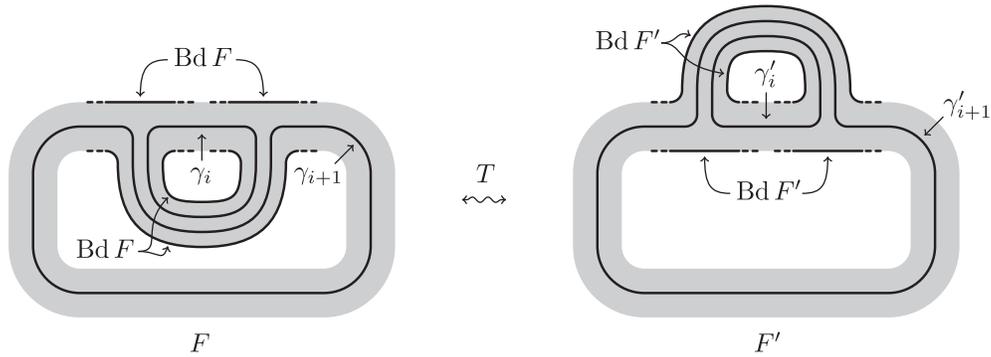

Figure 36. *The $T$ move.*

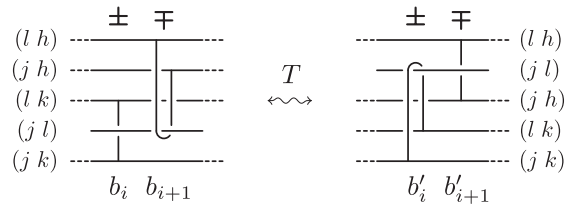

Figure 37. *A labelled line diagram presentation of the $T$ move.*

On the left, we have the regular fibre $F = F_0 \#_\eta F_t$ of $f$, with two consecutive twists $\gamma_i$ and $\gamma_{i+1}$ in a monodromy sequence for $f$. These twists have opposite sign and the corresponding vanishing cycles run parallel along the depicted band attached to the annulus $G_0$, which is not traversed by any other vanishing cycle of $f$. Then, the move consists in replacing $F$ by the regular fibre $F' = F_0 \#_{\eta'} F_t$ of $f'$, with the band attached on the opposite side of $G_0$, and the twists $\gamma_i$ and $\gamma_{i+1}$ with the twists $\gamma'_i$ and $\gamma'_{i+1}$ having the same signs. All the other twists in the monodromy sequence are left unchanged, but now they are thought of as twists in $F'$ instead of $F$ (this is possible, since the corresponding vanishing cycles are disjoint from the changed band). Hence, also the new band in $F'$ is traversed only by the vanishing cycles of $\gamma'_i$ and $\gamma'_{i+1}$.

It is not difficult to realize that the $T$ move preserves allowability of Lefschetz fibrations. In the allowable case, a labelled line diagram presentation of the $T$ move is depicted in Figure 37. Here, $j, l, k$, and $h$ are all different and only the (portions of the) horizontal lines involved in the move are drawn, while no vertical arc except those in the figure starts from or crosses behind the horizontal line labelled $(j\,l)$. The verification that the labelled line diagrams in the figure do actually represent the configurations in Figure 36 is straightforward and left to the reader.

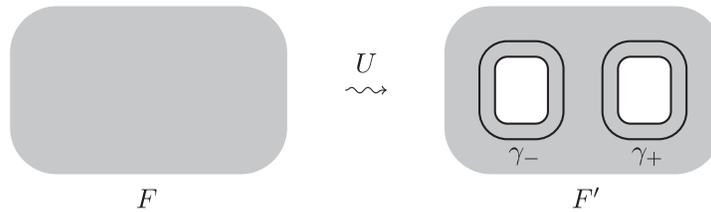

Figure 38. *The U move.*

Looking again at Figure 36, we see that the handlebodies $H_f$ and $H_{f'}$ are 2-equivalent. In fact, by sliding the 2-handle relative to $\gamma_{i+1}$ over that relative to $\gamma_i$ in $H_f$ and then cancelling the latter with the 1-handle generated by the band, we get the same handlebody resulting from the analogous operations performed on $H_{f'}$. Alternatively, it would be enough to observe that the diagrams in Figure 37 represent labelled ribbon surfaces related by labelled 1-isotopy, thanks to Proposition 6.

We remark that this argument still works even if the signs of the original twists $\gamma_i$ and $\gamma_{i+1}$ are not opposite. Moreover, it could be easily adapted to a generalized version of the $T$ move, where in place of the single twists $\gamma_{i+1}$ and $\gamma'_{i+1}$ there are sequences of twists $\gamma_{i+1}, \gamma_{i+2}, \ldots, \gamma_{i+k}$ in $F$ and $\gamma'_{i+1}, \gamma'_{i+2}, \ldots, \gamma'_{i+k}$ in $F'$. The only condition required for these twists is that the corresponding vanishing cycles run parallel to $\gamma_i$ and $\gamma'_i$ in the same order along the bands attached to $G_0$ in Figure 36, each cycle being allowed to traverse the bands more than once.

*U move.* This move will only be used to transform a non-allowable Lefschetz fibration into an allowable one, while it is not needed in the context of allowable Lefschetz fibrations. Given any Lefschetz fibration $f: W \to B^2$, it consists of making two holes in the interior of the regular fibre $F$ and then adding two singular points with vanishing cycles parallel to the new boundary components of $F$ and opposite signs, as shown in Figure 38.

More precisely, if $D \subset F$ is a disc disjoint from all the vanishing cycles of $f$, we put $F_0 = \mathrm{Cl}(F - D)$ and denote by $f_0: W_0 \to B^2$ the Lefschetz fibration with regular fibre $F_0$ and the same monodromy as $f$ (with the Dehn twists acting on $F_0$ instead of $F$). Then, we consider the fibre gluing $f' = f_0 \#_\eta t$, where the gluing map $\eta$ identifies a collar of $\mathrm{Bd}\, D$ in $F_0$ with a collar of the outer boundary component of $F_t$ in Figure 34.

We call a $U$ move the modification $U: f \rightsquigarrow f' = f_0 \#_\eta t$, as well as its inverse $U^{-1}: f' \rightsquigarrow f$. We notice that, the Lefschetz fibration $f'$ has the desired regular fibre $F'$, while a mapping monodromy sequence for it can be obtained from one for $f$ by inserting anywhere in the sequence a pair $\gamma_+$ and $\gamma_-$ of a positive and a negative twist parallel to the new boundary components of $F'$ (see Figure 38).

The handlebody $H_{f'}$ can be shown to be 2-equivalent to $H_f$, by applying to it the same 2-handle sliding and handle cancellation considered above in the case of the $T$ move.

If $f$ is a non-allowable Lefschetz fibration with regular fibre $F$, then we can perform on it a sequence of $U$ moves, to make $F$ into a bounded surface and/or to insert a hole in the interior of each subsurface of $F$ bounded by a vanishing cycle, thus obtaining an allowable Lefschetz fibration.

We conclude this section by briefly discussing the independence of the moves defined above. First of all, we observe that the $U$ move is clearly independent from the others, being the only one that does not preserve allowability. Nevertheless, it can be generated by moves $S$ and $T$ in the context of allowable Lefschetz fibrations, as will follow from the results of Section 8. Therefore, the $U$ move should be considered as an auxiliary move, used just to get allowability.

The independence of the $S$ move from the $T$ move is easy to see. In fact, only the $S$ move changes the parity of the number of boundary components of the regular fibre. To establish that the $T$ move is independent from the $S$ move, we need to introduce the Euler class of a Lefschetz fibration and study how it is affected by moves.

Let $f: W \to B^2$ be a Lefschetz fibration. Then, the restriction $f_{|W'}: W' \to B^2$, with $W' = W - \{w_1, w_2, \ldots, w_n\}$ the complement of the singular set of $f$, is a submersion. Hence, we can consider the distribution $\xi_f$ of oriented planes on $W'$ given by the kernel of the tangent map $Tf_{|W'}: TW' \to TB^2$, and its Euler class $e(\xi_f) \in H^2(W')$. The *Euler class* of the Lefschetz fibration $f$ is defined as $e(f) = (i^*)^{-1}(e(\xi_f)) \in H^2(W)$, where $i^*: H^2(W) \to H^2(W')$ is the isomorphism induced by the inclusion $i: W' \subset W$.

Note that, since $TW' \cong \xi_f \oplus \xi_f^\perp$ with $\xi_f^\perp \cong f_{|W'}^*(TB^2)$ a trivial bundle, the mod 2 reduction of $e(f)$ coincides with $w_2(W)$, the second Stiefel–Whitney class of $W$.

Now, we want to express the Euler class $e(f)$ in terms of a mapping monodromy sequence $(\gamma_1, \gamma_2, \ldots, \gamma_n)$ for $f$, when the regular fibre $F \cong F_{g,b}$ of $f$ is a bounded surface, that is, $b > 0$. In this case, $TF_{g,b}$ is trivial and we can choose a positive frame field $(u_1, u_2)$ on $F_{g,b}$. Moreover, a basis for the cellular 2-chain group $C_2(W)$ is provided by the cores of the 2-handles of the handlebody decomposition $H_f$ of $W$ with any given orientation, whose boundaries are the vanishing cycles $c_1, c_2, \ldots, c_n \subset F_{g,b}$ with the induced orientation. We use the same notation $c_1, c_2, \ldots, c_n$ for those generators of $C_2(W)$ and denote by $c_1^*, c_2^*, \ldots, c_n^*$ the dual generators of the cellular 2-cochain group $C^2(W)$.

PROPOSITION 12. *Given a Lefschetz fibration $f: W \to B^2$ with bounded regular fibre $F \cong F_{g,b}$ and oriented vanishing cycles $c_1, c_2, \ldots, c_n \subset F_{g,b}$, and any positive frame field $(u_1, u_2)$ for $F_{g,b}$, we have $e(f) = [\varepsilon]$ with $\varepsilon = \sum_{i=1}^n \mathrm{rot}(c_i) c_i^* \in C^2(W)$, where $\mathrm{rot}(c_i)$ is the rotation number of $c_i$ with respect to $(u_1, u_2)$.*

*Proof.* We start with a handlebody decomposition $H_f = W_1 \cup H_1^2 \cup H_2^2 \cup \cdots \cup H_n^2$ of $W$, where $W_1$ is a 1-handlebody decomposition of $B^2 \times F_{g,b}$ (cf. proof of Proposition 9 and take into account that $\mathrm{Bd}\, F \neq \emptyset$). Then, each 2-handle $H_i^2$ contains one singular point $w_i$ and is modelled on the Hopf fibration $h_\pm: H_\pm \to B^2$. Hence, the restriction to $H_i^2 - \{w_i\}$ of the plane field $\xi_f$ is the pull-back of the plane field $\xi_{h_\pm}$ on $H_\pm - \{0\}$ under the fibred equivalence $H_i^2 \cong H_\pm$.

Recalling that $H_\pm = \{(z_1, z_2) \in \mathbb{C} \mid |z_1|^2 + |z_2|^2 \leqslant 2 \text{ and } |h_\pm(z_1, z_2)| \leqslant 1\}$, with $h_+(z_1, z_2) = z_1^2 + z_2^2$ and $h_-(z_1, z_2) = z_1^2 + \bar{z}_2^2$ for all $(z_1, z_2) \in \mathbb{C}^2$, and putting $z_1 = x_1 + iy_1$ and $z_2 = x_2 + iy_2$, a straightforward computation shows the following.

(1) The vanishing cycle $c$ of $h_\pm$ in the regular fibre $F_1 = h_\pm^{-1}(1)$ is (up to isotopy) the circle of equations $x_1^2 + x_2^2 = 1$ and $y_1 = y_2 = 0$;
(2) A trivializing positive frame field $(v_1, v_2)$ for $\xi_{h_\pm}$ is given by

$$v_1 = -x_2 \frac{\partial}{\partial x_1} \mp y_2 \frac{\partial}{\partial y_1} + x_1 \frac{\partial}{\partial x_2} \pm y_1 \frac{\partial}{\partial y_2} \quad \text{and} \quad v_2 = y_2 \frac{\partial}{\partial x_1} \mp x_2 \frac{\partial}{\partial y_1} \mp y_1 \frac{\partial}{\partial x_2} + x_1 \frac{\partial}{\partial y_2};$$

(3) The restriction of $v_1$ to $c$ is a positive tangent vector field on $c$ with the usual counterclockwise orientation.

On the other hand, on $W_1$ we consider the trivialization of $\xi_f$ induced by the pull-back of the frame field $(u_1, u_2)$ under the projection $\pi: W_1 \cong B^2 \times F_{g,b} \to F_{g,b}$, which we still denote by $(u_1, u_2)$. Property 3 of the frame field $(v_1, v_2)$ implies that, for every $i = 1, 2, \ldots, n$, the rotation number $\mathrm{rot}(c_i)$ represents the obstruction to extending the frame field $(u_1, u_2)$ over $H_i^2 - \{w_i\}$, since this comes from the fibre gluing of $H_\pm - \{0\}$ to $W_1$ along $c_i$ (see proof of Proposition 9).

Thus, the cohomology class of $\varepsilon$ in $H^2(W')$ coincides with $e(\xi_f)$, and the proposition follows at once. □

Proposition 12 allows us to easily compute the changes in the Euler class $e(f)$ induced by any move performed on $f$. This is done in the following proposition, which obviously implies the independence of the $T$ move from the $S$ move.

PROPOSITION 13. *Let $f: W \to B^2$ and $f': W \to B^2$ be Lefschetz fibrations with bounded regular fibres. If $f'$ is obtained from $f$ by an $S$ move, then $e(f') = e(f)$. While, if $f'$ is obtained from $f$ by a $T$ move or a $U$ move, then $e(f') = e(f) + 2[\delta]$ for some (generically cohomologically non-trivial) cocycle $\delta \in C^2(W)$.*

*Proof.* We choose trivializing frame fields $(u_1, u_2)$ and $(u_1', u_2')$ for $F$ and $F'$, respectively, to coincide with the standard one in Figures 33, 36 and 38, and assume all the vanishing cycles in those figures to be counterclockwise oriented. Then we use Proposition 12 to evaluate the difference $e(f') - e(f)$.

In the case when $f'$ is a Hopf-stabilization of $f$, we have that $e(f') - e(f) = [c^*]$, where $c^*$ is the dual of the generator $c$ of $C_2(W)$ corresponding to the new Dehn twist $\gamma_\pm$ in Figure 33. But this is cohomologically trivial in $W$, hence $e(f') = e(f)$.

If $f$ and $f'$ are related by a $T$ move as in Figure 36, then $e(f') - e(f) = [(c_i')^* + (c_{i+1}')^*] - [c_i^* + c_{i+1}^*]$. Since $c_i^*$ and $(c_i')^*$ are, respectively, cohomologous to $c_{i+1}^*$ and $-(c_{i+1}')^*$ in $W$, we can write $e(f') - e(f) = 2[\delta]$ with $\delta = -c_i^*$. Similarly, if $f$ and $f'$ are related by a $U$ move as in Figure 38, we have $e(f') - e(f) = 2[\delta]$ with $\delta = c_+^*$. □

REMARK 14. *Any integral lifting of the second Stiefel–Whitney class of a four-dimensional 2-handlebody $W$ is the Euler class of a Lefschetz fibration $W \to B^2$.*

Since any two such liftings differ by an even class (by the universal coefficient theorem), that statement follows once we show that starting from any Lefschetz fibration $f: W \to B^2$ with bounded fibre $F$, one can construct another Lefschetz fibration $f': W \to B^2$ such that $e(f') = e(f) \pm 2c^*$ for an arbitrary vanishing cycle $c \subset F$ of $f$. Up to isotopy, we can assume that there is a disc in Int $F - c$ and an arc joining it to Bd $F$, which intersects $c$ transversely in a single point and do not intersect any other vanishing cycles. Then by performing a $U$ move inside this disc, we get $e(f') = e(f) - 2c^*$. Changing $f'$ by two more $U$ moves, performed in a similar way with the two new vanishing cycles introduced by the former $U$ move in place of $c$, we get $e(f') = e(f) + 2c^*$.

## 8. *The main theorem*

As a preliminary step to prove the equivalence theorem for Lefschetz fibrations, we need to translate into the language of rectangular diagrams the moves for labelled flat planar diagrams considered in Proposition 8. In doing that, we continue to use for rotated moves the 'prime notation' and the graphical 'rounded bottom-left corner' rule introduced in Section 4.

First of all, we consider the plane isotopy moves in Figure 39. These are intended to relate planar rectangular diagrams, which are isotopic in the projection plane (crossing and ribbon intersections are assumed to be preserved by the isotopy). We observe that only two rotations of the moves $r_1$, $r_3$ and $r_6$ are enough, due to their symmetry. Move $r_1$ switches two horizontal bands which are contiguous in the vertical order, under the assumption that the horizontal intervals they span do not overlap. Similarly, move $r_1'$ switches two contiguous vertical bands under the analogous assumption. Here we include, as degenerate horizontal or vertical bands, the end points of vertical or horizontal ones, respectively. In all the moves, when a band coming out from the box is translated, we assume that the translation is small enough to not interfere with the rest of the diagram.

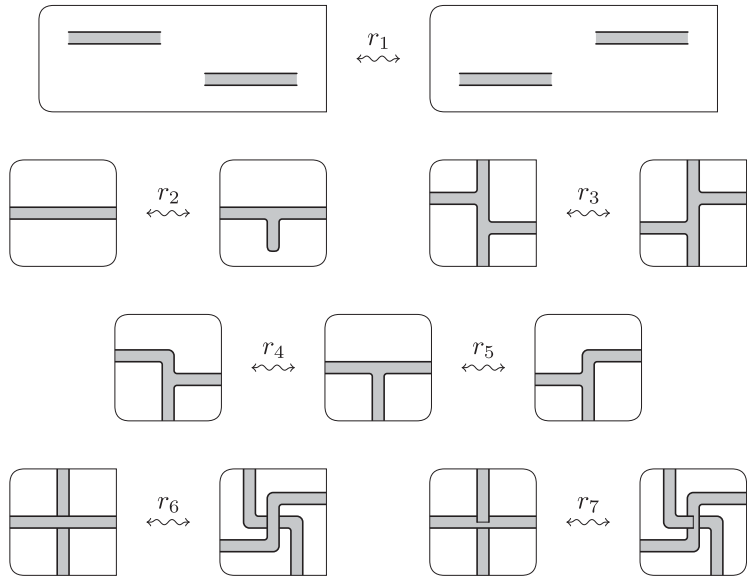

Figure 39. *Plane isotopy moves for rectangular diagrams.*

Figures 40–42 provide the rectangular versions of the moves in Figures 3, 6, 7 and 24. We observe that, no rotated move is needed here, thanks to the moves in Figure 39 that allow us to rotate by $k\pi/2$ any local configuration in the diagram. Moreover, we put each move in the most convenient form for applying Rudolph's braiding procedure to it.

Lemma 15. *Two labelled rectangular diagrams represent 2-equivalent connected four-dimensional 2-handlebodies as simple branched coverings of $B^4$ if and only if they are related by rectangular (de)stabilization and the moves $r_1$ to $r_{25}$ in Figures 39–42.*

*Proof.* Moves $r_1$ to $r_7$ in Figure 39, together with their allowed rotated versions, suffice to realize any plane isotopy between rectangular diagrams. In fact, any flat planar diagram is uniquely determined by its planar core graph, together with some extra information on the vertices corresponding to crossings and ribbon intersections. Now, it is not difficult to realize that any isotopy of the graph, as well as any deformation of its structure, can be approximated by using the moves in Figure 39. In particular, the two moves $r_{26}$ and $r_{27}$ in Figure 43, which are obviously needed in order to approximate isotopies along an edge, can be obtained from $r_4$ and $r_5$ modulo $r_2$.

Then, it suffices to show that the remaining moves $r_8$ to $r_{25}$ generate all the moves listed in Proposition 8, in the presence of moves $r_1$ to $r_7$.

The moves in Figure 40 together with $r_2$ and $r_3$, generate the moves $s_5$ to $s_{18}$ in Figures 6 and 7 and the move $s_{27}$ in Figure 27. The only non-trivial facts in this respect are the following: (1) move $s_9$, which does not have an explicit representation in Figure 40, can be obtained as a composition of moves $r_2$, $r_8$ and $r_{10}$; (2) $r_{14}$ and $r_{15}$ completed with the right terminations, give $s_{15}$ and $s_{16}$ (after contracting the extra tongue, by using $r_2$, $r_4$, $r_5$, $r_8$ and $r_9$); (3) similarly, $r_{16}$, $r_{17}$, $r_{18}$ and $r_{19}$ give $s_5$, $s_6$, $s_{17}$ and $s_{18}$.

While the moves $r_{20}$ to $r_{23}$ in Figure 41 and the moves $r_{24}$ and $r_{25}$ in Figure 42 are rectangular versions of the 1-isotopy moves $s_1$ to $s_4$ in Figure 3 and of the covering moves $c_1$ and $c_2$ in Figure 27, respectively. □

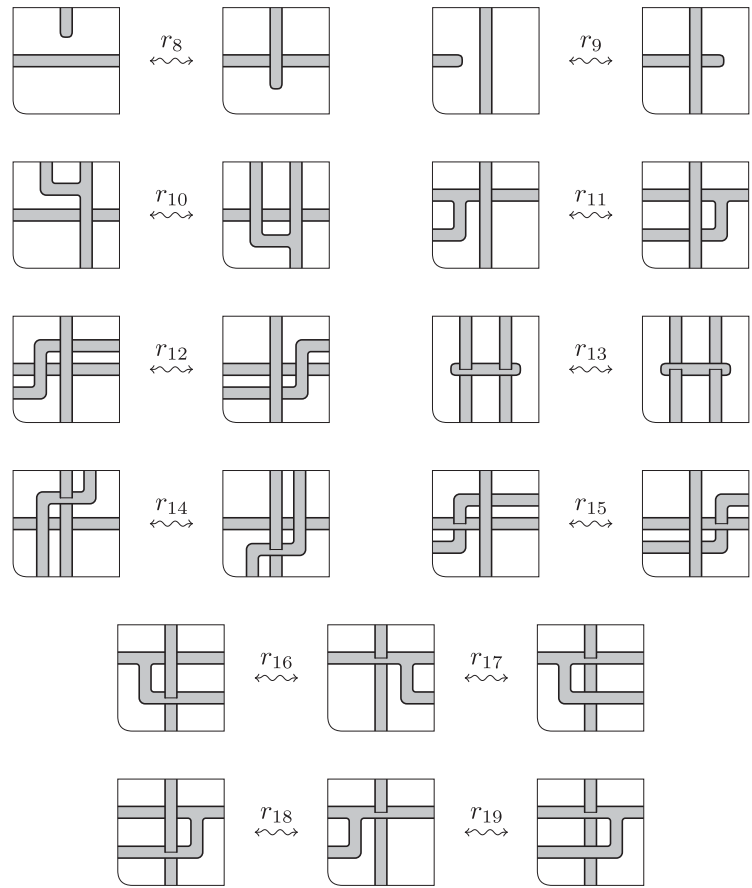

Figure 40. *Three-dimensional isotopy moves for rectangular diagrams.*

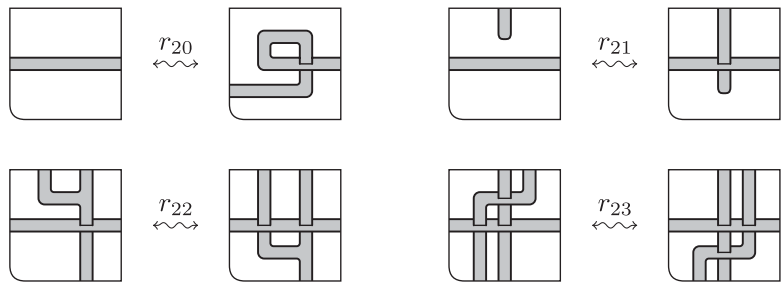

Figure 41. *1-isotopy moves for rectangular diagrams.*

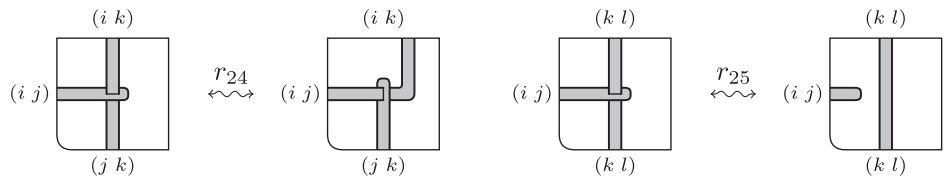

Figure 42. *Covering moves for colored rectangular diagrams.*

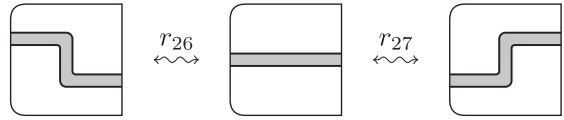

Figure 43. *Breaking edges.*

At this point, we are in a position to prove our main theorem.

Theorem A.  *Any two allowable Lefschetz fibrations $f: W \to B^2$ and $f': W' \to B^2$ represent 2-equivalent four-dimensional 2-handlebodies $H_f$ and $H_{f'}$ if and only if they are related by fibred equivalence and the moves $S$ and $T$. Moreover, the allowability hypothesis can be relaxed by using in addition move $U$.*

*Proof.* As we have seen in Section 7, fibred equivalence and moves $S$, $T$ and $U$ do not change the 2-equivalence class of the total space of a Lefschetz fibration. This gives the 'if' part of the statement. Concerning the 'only if' part, the reduction to the allowable case immediately follows from the fact that any Lefschetz fibration can be made into an allowable one by performing $U$ moves on it (see Section 7).

Now, let $f$ and $f'$ be allowable Lefschetz fibrations as in the statement. Thanks to Proposition 11 and the consideration following it, they can be represented by band presentations $S$ and $S'$ of certain labelled braided surfaces. Moreover, by Proposition 3 (cf. discussion on $S$ move in Section 7), up to positive stabilization and fibred equivalence of Lefschetz fibrations, the labelled surfaces $S$ and $S'$ can be assumed to have monotonic bands. Finally, we perform on these surfaces the flattening procedure described at the end of Section 3, getting in this way two labelled rectangular diagrams, which we still denote by $S$ and $S'$. Lemma 15 tells us that $S$ and $S'$ are related by rectangular (de)stabilization and moves $r_1$ to $r_{25}$.

Concerning rectangular (de)stabilization, we observe that it can be realized by the addition/deletion of a separate short horizontal band labelled $(i\, d+1)$ (cf. Section 5). But, once the braiding procedure has been applied, this band results into a separate sheet labelled $(i\, d+1)$, which affects the covering but not the Lefschetz fibration (cf. discussion on $S$ move in Section 7).

The rest of the proof is devoted to showing that Rudolph's braiding procedure makes moves $r_1$ to $r_{25}$ (including the rotated versions of moves $r_1$ to $r_7$) into modifications of labelled braided surfaces, which can be realized by labelled band sliding (that means labelled braided isotopy), (de)stabilization of labelled braided surfaces (which corresponds to the $S$ move) and the labelled braided surface representation of the $T$ move (see Figure 37). Note that the braiding procedure applied to the rectangular diagrams $S$ and $S'$ gives back the original braided surfaces (cf. Proposition 4).

We will include in the argument also the moves $r_{26}$ and $r_{27}$ (and their rotated versions). These are not strictly needed, but they help us to simplify the handling of the other moves. Moves will be considered in a convenient order not consistent with the numbering.

In the computations below, we will use the line notation for the rectangular diagrams, as we already did for the braided surfaces. This consists of a rectangular diagram of the core graph of the represented ribbon surface, with the conventions described in Figure 44 for the terminal bands and the ribbon intersections.

*Moves $r_1, r_{26}$ and $r_{27}$.* We observe that the two bands of the rectangular diagram involved in these moves give rise to two sets of bands in the braided surface, due to the preliminary replacements of local configurations that are not in the restricted form for Rudolph's procedure.

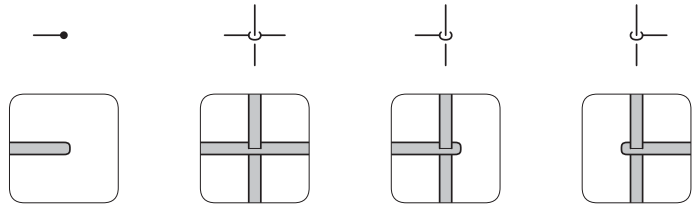

Figure 44. *The line notation for rectangular diagrams.*

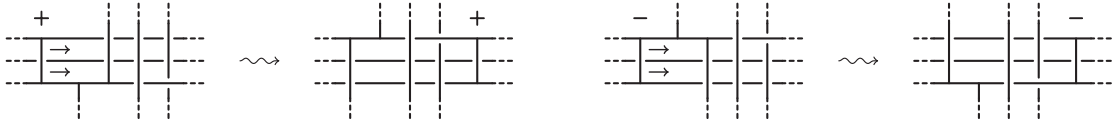

Figure 45. *Sliding a band to the right.*

Actually, those replacements make a single move $r_1$ between the original rectangular diagrams into a finite sequence of moves $r_1$ between the corresponding diagrams in restricted form. Analogously, a single move $r_{26}$ or $r_{27}$ is made into a sequence of moves $r_1$ and one move $r_{26}$ or $r_{27}$, respectively. Therefore, we can limit ourselves to consider moves $r_1, r_{26}$ and $r_{27}$ between diagrams in restricted form, each horizontal band of which generates a single sheet of the braided surface.

In the case of move $r_1$, from the two bands on the left we get two sheets of the braided surface, such that the upper or the lower one is trivial, respectively, on the right or on the left of a certain abscissa. We stabilize the braided surface by inserting a trivial sheet immediately over those two sheets and connecting it to the lower one by a positive half-twisted band located at that abscissa. After that, we make the lower sheet trivial by sliding this band to the right and then we remove it by destabilizing. The left side of Figure 45 shows the effect of the sliding on the four possible types of bands we can meet. The resulting braided surface is the one given by the diagram on the right side of move $r_1$.

Move $r_{26}$ can be treated in a similar way. Here, we just slide to the right the positive half-twisted band of the braided surface corresponding to the vertical band of the diagram on the left side of the move and then remove the lower sheet by destabilizing. The left side of Figure 45 still describes the sliding in this case if we ignore the sheet in the middle.

Move $r_{27}$ is the up-down symmetric of $r_{26}$. Hence, the same argument holds for it, except for the half-twisted band being negative and the sliding of it working as in the right side of Figure 45.

*Moves $r'_1, r'_{26}$ and $r'_{27}$.* The reduction to the case of diagrams in restricted form goes as above (rotate everything in the reasoning). In this case, move $r'_1$ just means a band sliding (we are changing the order of commuting bands), while moves $r'_{26}$ and $r'_{27}$ can be interpreted as a sliding followed by a destabilization, as in Figure 46.

From now on, we will implicitly use the above moves to localize the other ones, by breaking all the bands coming out from the involved local configurations. In this way, we can disregard the small translations of those bands needed for the move to take place. Moreover, the replacements needed to get configurations in restricted form can be performed in any order (cf. Section 4) and the moves can be assumed to directly act on diagrams in restricted form, after the replacements.

*Moves $r_2, r'_2, r''_2, r'''_2$.* Figure 47 shows how to deal with these moves. For $r_2$, we have only a destabilization, for $r'_2$ and $r'''_2$ we first need to perform one sliding, while for $r''_2$ we can destabilize the top sheet and then continue as for $r_{26}$.

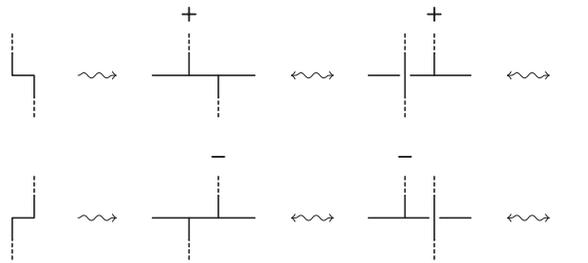

FIGURE 46. Moves $r'_{26}$ and $r'_{27}$.

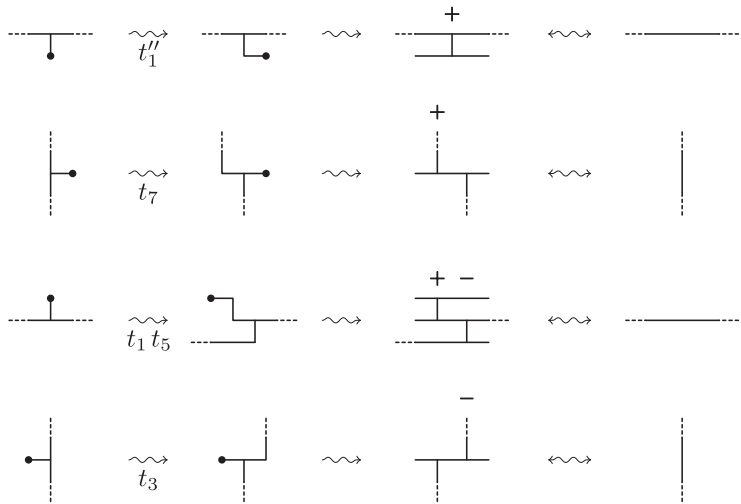

FIGURE 47. Moves $r_2, r'_2, r''_2$ and $r'''_2$.

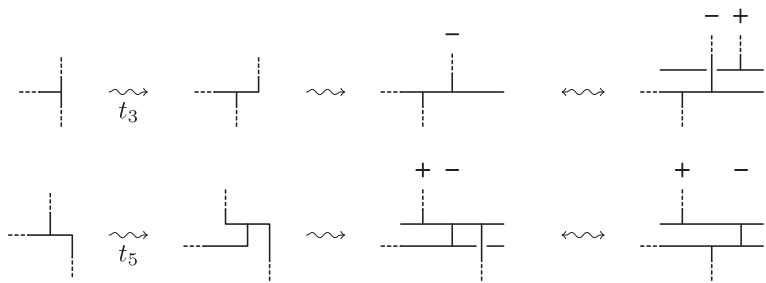

FIGURE 48. Move $r'''_5$.

*Moves* $r_4, r'_4, r''_4, r'''_4$ *and* $r_5, r'_5, r''_5, r'''_5$. After the replacements in Figure 23: move $r_4$ follows from $r_{26}$; move $r'_4$ follows from $r'_1$ and $r'_{26}$; move $r_5$ follows from $r_{27}$; moves $r'''_4, r'_5$ and $r''_5$ are tautological. Moreover, moves $r''_4$ and $r'''_5$ are equivalent modulo $r_{26}$. Finally, move $r'''_5$ is considered in Figure 48. Here, a stabilization and a band sliding are performed at the end of the top and the bottom line respectively, and the two resulting rightmost diagrams are equivalent up to band sliding.

*Moves* $r_3, r'_3$. These moves can be easily seen to be equivalent modulo moves $r_4, r_5, r_{26}, r_{27}$ and their rotated versions. Move $r'_3$ is considered in Figure 49. Here, the two rightmost diagrams are equivalent up to band sliding.

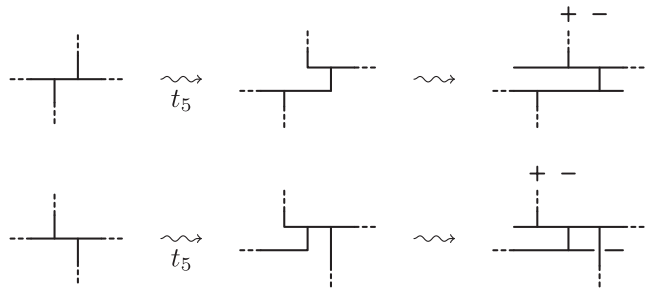

FIGURE 49. *Move $r_3'$.*

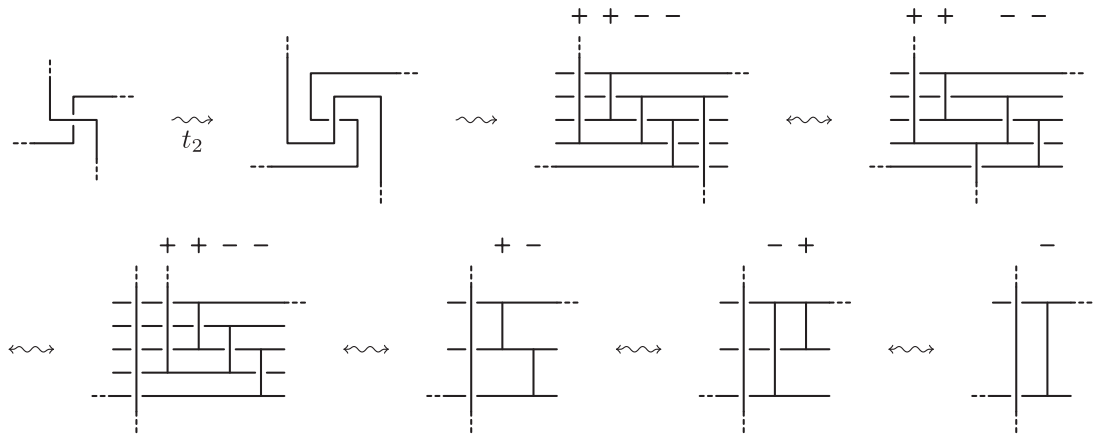

FIGURE 50. *Move $r_6'$.*

*Moves $r_6, r_6'$.* We observe that move $r_6$ is tautological, since it coincides with $t_2$. Move $r_6'$ is treated in Figure 50. Also in this case band sliding and (de)stabilization suffice. In particular, after the last step in the figure, the band between the two remaining sheets has to be slid up to the right end (as in Figure 45), to allow a final destabilization.

*Moves $r_8, r_9$ and $r_{21}$.* For $r_8$ and $r_{21}$ it suffices to note that, once the replacement $t_1''$ is applied to the terminal vertical band, the sheet deriving from this band can be removed by destabilization on both the sides of the move. Move $r_9$ does not affect at all the resulting braided surface.

*Moves $r_{11}, r_{12}, r_{15}, r_{16}, r_{18}$ and $r_{20}$.* The local configurations involved in all these moves are in the restricted form (hence, no replacement is required) and it is easy to see that they just change the braided surface by a sliding of one of the two half-twisted bands which appear over the other one. In particular, these bands commute for the first three moves. Actually, for move $r_{20}$ some further band sliding and a destabilization are required, like for $r_{26}$ and $r_2''$.

*Moves $r_{17}$ and $r_{19}$.* In Figure 51, we compare the two braided surfaces originated by the local configurations in move $r_{17}$. Once again we see that they are equivalent up to band slidings and stabilization.

The case of move $r_{19}$ is symmetric to that of $r_{17}$ modulo moves $r_1, r_1', r_{26}, r_{26}'$ and $r_6'$. In fact, in Figure 52, we show how the right side of move $r_{19}$ can be put in restricted form by using those moves (instead of the replacement $t_2$). The result is symmetric to the second diagram in Figure 51.

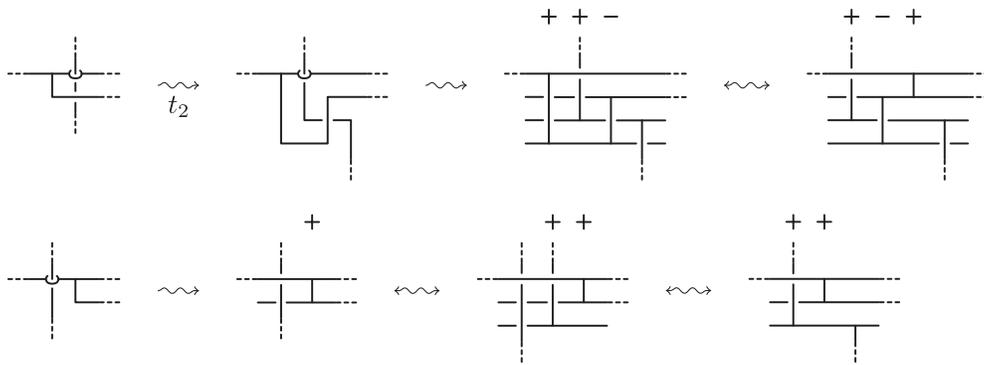

FIGURE 51. *Move $r_{17}$.*

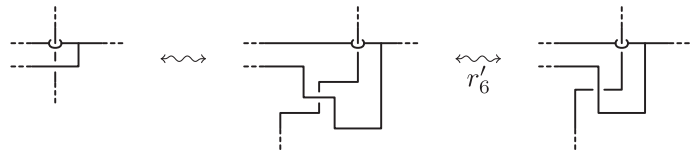

FIGURE 52. *Move $r_{19}$.*

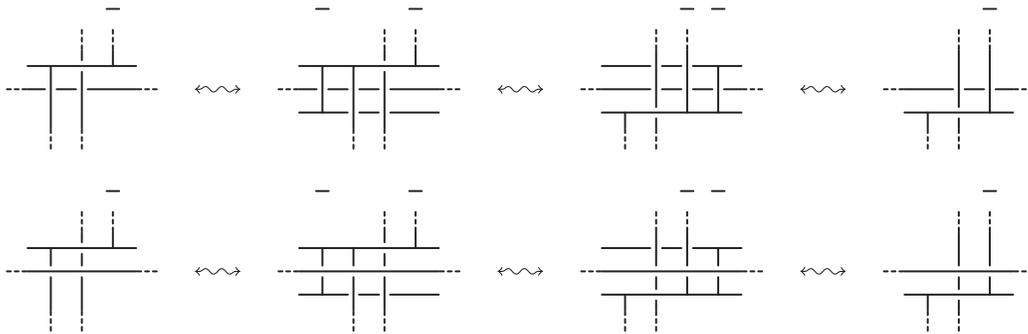

FIGURE 53. *Moves $r_{14}$ and $r_{23}$.*

*Moves $r_{10}, r_{14}, r_{22}$ and $r_{23}$.* The argument for all these moves is essentially the same. Moves $r_{14}$ and $r_{23}$ are already in restricted form, so no replacement is needed. For the moves $r_{10}$ and $r_{22}$, the same replacement $t_3$ occurs on both sides. In any case, once the move is in restricted form, we have to change the position of the sheet corresponding to the short horizontal band from top to bottom. Up to stabilization, this can be done by band sliding, as shown in Figure 53 for moves $r_{14}$ and $t_{23}$. The procedure for the other moves is analogous.

Before passing to the remaining moves, we introduce the auxiliary 1-isotopy moves for rectangular diagrams depicted in Figure 54. Here, it does not matter what the labelling is.

Moves $r_{28}, r_{29}$ and $r_{30}$ are nothing else than rectangular versions of different planar projections of moves $s_2$, $s_3$ and $s_4$, respectively, hence we know that they follow from the moves $r_1$ to $r_{23}$ (in particular move $r_{13}$ is needed here). But, since these auxiliary moves will be used to deal with moves $r'_7$ and $r_{13}$, we directly consider the effect on the braided surface resulting from the braiding procedure applied to them, independently of the other moves. Actually, the reader can easily check that the modifications they induce on the corresponding braided surface are completely analogous to those induced by moves $r_{21}, r_{22}$ and $r_{23}$, respectively, and can be realized by (de)stabilization and band sliding as well.

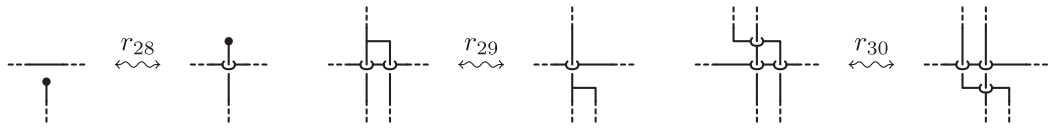

Figure 54. *Auxiliary moves.*

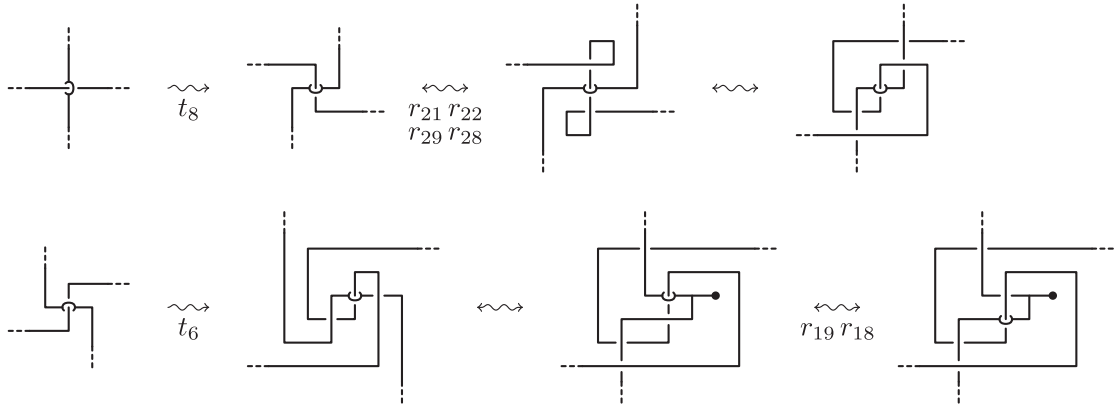

Figure 55. *Move $r'_7$.*

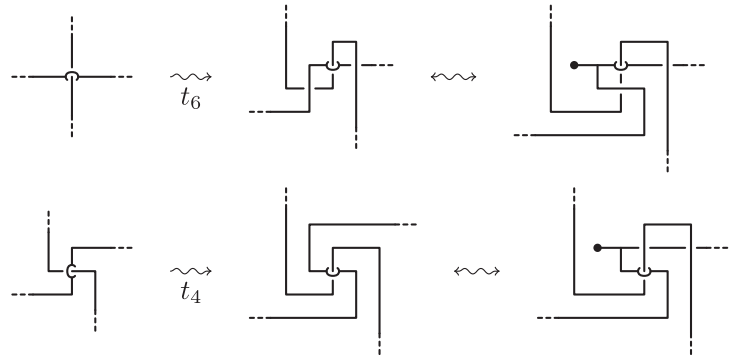

Figure 56. *Move $r''_7$.*

*Moves $r_7, r'_7, r''_7, r'''_7$.* Concerning $r_7$, we see that once the replacement $t_8$ is applied to the diagram on the right side of the move, it turns out to be equivalent to that on the left side up to moves $r_1, r'_1, r_{27}, r'_{27}$. Move $r'''_7$ is tautological, since it coincides with $t_4$. So, we are left with moves $r'_7$ and $r''_7$. These are, respectively, treated in Figures 55 and 56, where they are derived from the moves considered above after the required replacements. In particular, the third diagram in the first line of Figure 55 can be proved to be equivalent to the second one, by cancelling the two opposite kinks, once the top one has been moved down passing through the horizontal band in the middle. The last operation can be realized in a straightforward way, by using the moves already considered above and the auxiliary moves $r_{28}$ and $r_{29}$.

*Moves $r_{24}$ and $r_{25}$.* We consider these moves in Figures 57 and 58. Here, we start with one side of the move (the right side for $r_{24}$ and the left side for $r_{25}$) and end up with the braided surface corresponding to the other side. In both the figures, all the modification of the labelled braided surface are (de)stabilizations or band slidings, except the first step in the second line. This step consists of applying three negative half-twists on the interval between the two sheets to the band coming from the top in Figure 57, and two positive half-twists on the interval between the two sheets to the band coming from the bottom in Figure 58. Since both the labelled braids corresponding to those multiple half-twists belong to the kernel of the lifting homomorphism (cf. Section 5 and [**17**]), up to fibred equivalence the represented Lefschetz fibration does not change.

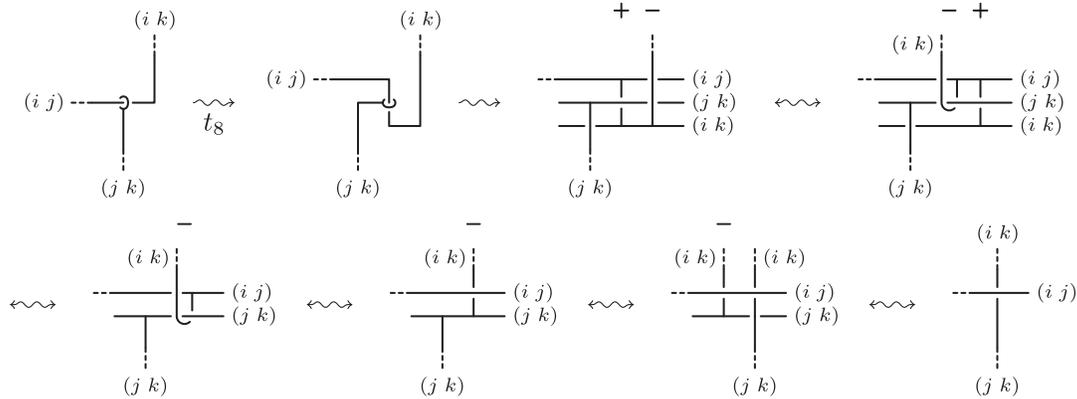

Figure 57. *Move $r_{24}$.*

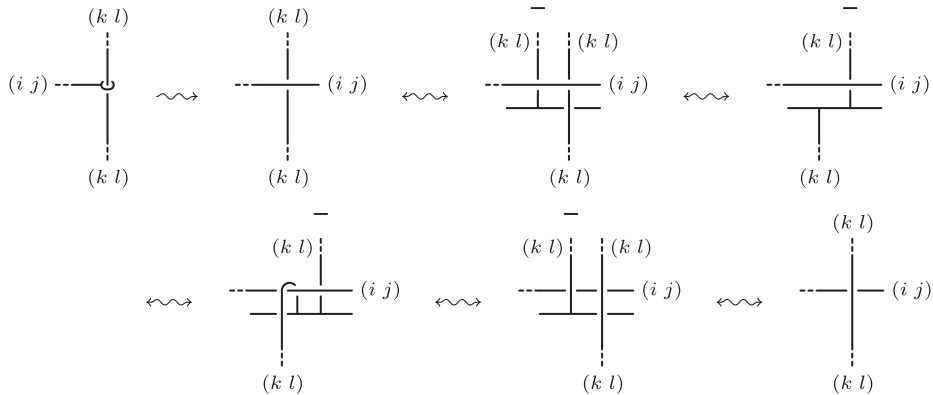

Figure 58. *Move $r_{25}$.*

*Move $r_{13}$.* First of all, we show that modulo the other moves, move $r_{13}$ can be reduced to the case when the labelling satisfies very restrictive conditions. Namely, if $(i\ j)$ is the label of the horizontal disc and $\sigma$ and $\tau$ are the bottom labels of the vertical bands passing through it, then we can assume $\sigma = (i\ k)$ and $\tau = (j\ l)$.

To see this, we first consider the rectangular move $r_{31}$ in Figure 59. We observe that, replacing back the annuli with the corresponding half-twists, this move just slide the ribbon intersection formed by the vertical band with the horizontal one, from right to left along the latter across a single half-twist. Then, this move is in fact a three-dimensional diagram move and it holds for any labelling of the vertical and horizontal bands. But deriving $r_{31}$ in this general form would involve move $r_{13}$, while with the labelling specified in Figure 59 it can be derived without using $r_{13}$ and still suffices for our purposes.

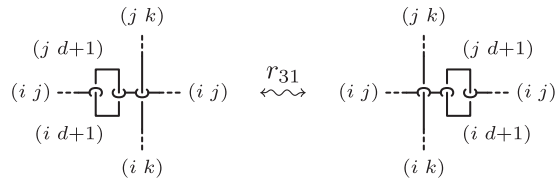

FIGURE 59. *The auxiliary move* $r_{31}$.

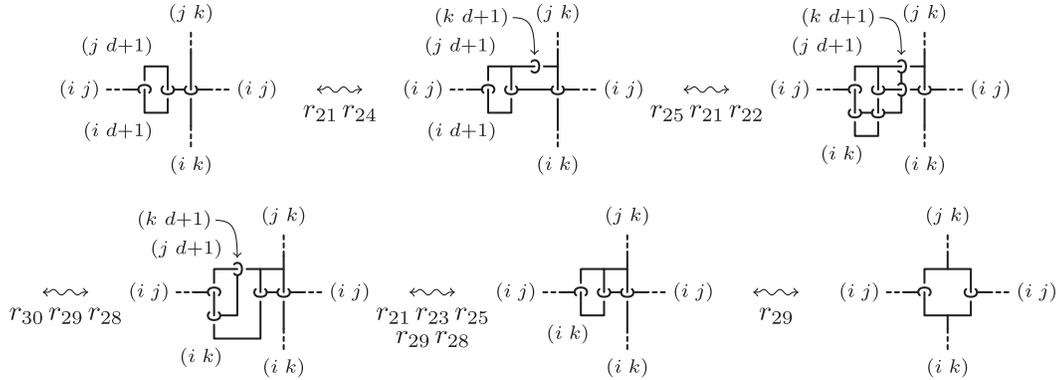

FIGURE 60. *Deriving the auxiliary move* $r_{31}$.

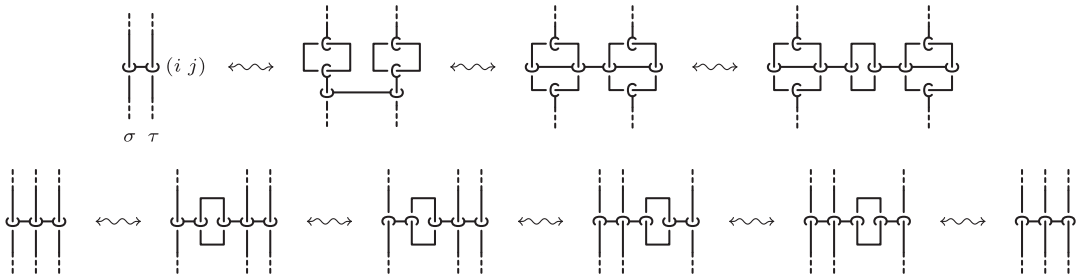

FIGURE 61. *Reducing move* $r_{13}$ *to the case* $\sigma = (i\ k)$ *and* $\tau = (jl)$.

Figure 60 shows how to derive move $r_{31}$ with that labelling from the other moves except $r_{13}$. Here, the moves indicated under the arrows are intended up to the planar isotopy moves in Figure 39 (in particular up to rotations).

At this point, we are in a position to deal with move $r_{13}$. We begin by modifying both sides of the move as in the top line of Figure 61, where the disc and the bands are broken by using move $c_4$ in Figure 25 and some 1-isotopy moves are performed (including $r_{28}$ and $r_{30}$). In this way, the original move is changed into two reversing moves involving three ribbon intersections, but now the bottom labels of the vertical bands are all different from each other and from $(i\ j)$. In the second line of Figure 61, we realize such a reversing move in terms of two moves $r_{13}$ with the same constraints on labelling. This also requires two moves $c_4$ at the first and last step and a move $r_{31}$ at the middle step.

By this argument, we can assume that also $\sigma$ and $\tau$ in the original move $r_{13}$ are different from each other and from $(i\ j)$. Under such assumption, if $\sigma$ or $\tau$ are disjoint from $(i\ j)$, then the move can easily be deduced from the covering moves in Figure 42. In the remaining cases, up to symmetry we have $\sigma = (i\ k)$ while $\tau$ can be one of $(i\ l)$, $(j\ k)$ or $(j\ l)$ (here, we assume $i, j, k$ and $l$ all different). At this point, the reduction to the only case when $\tau = (j\ l)$ is immediate.

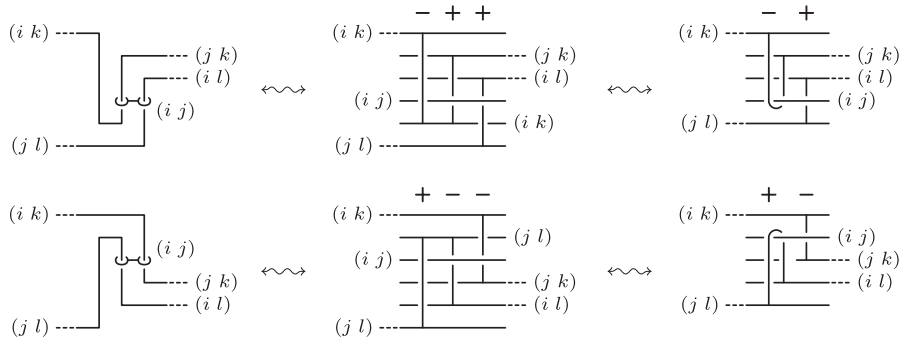

Figure 62. *Move $r_{13}$.*

Figure 62 shows how to realize move $r_{13}$ with $\sigma = (i\ k)$ and $\tau = (j\ l)$ in terms of the $T$ move described in Section 7. This is needed to relate the two rightmost braided surfaces, once they are obtained by suitable stabilizations (the sign of the stabilizing band is irrelevant).  □

Remark 16. According to Theorem A, any $U$ move between allowable Lefschetz fibrations can be generated by $S$ and $T$ moves. Actually, this could be proved directly by induction on the number of vanishing cycles that separate the region where the $U$ move is performed from the boundary of the fibre. Indeed, this number can be reduced to zero by suitable $T$ moves and slidings, and after that the $U$ move can be trivially realized by two opposite $S$ moves.

## 9. *Open books*

Given a closed connected oriented three-manifold $M$, by an *open book* structure on $M$, we mean a smooth map $f: M \to B^2$ such that the following properties hold.

(1) The restriction $f_{|T}: T = \mathrm{Cl}(f^{-1}(\mathrm{Int}\, B^2)) \to B^2$ is a (trivial) fibre bundle, with $L = f^{-1}(0) \subset M$ a smooth link, called the *binding* of the open book, and $T \subset M$ a tubular neighbourhood of $L$;

(2) The composition $\rho_f = \rho \circ f_{|M-L}: M - L \to S^1$, with $\rho: B^2 - \{0\} \to S^1$ the projection defined by $\rho(x) = x/\|x\|$, is a locally trivial fibre bundle, whose fibre is the interior $\mathrm{Int}\, F$ of a compact connected orientable bounded surface $F$, called the *page* of the open book.

Such a map induces an *open book decomposition* of $M$ into compact connected orientable bounded surfaces $F_s = f^{-1}([0, s]) = \mathrm{Cl}(\rho_f^{-1}(s)) \subset M$ with $s \in S^1$, called the *pages* of the decomposition. These are all diffeomorphic to $F$ and only meet at their common boundary $\mathrm{Bd}\, F_s = L$. Moreover, $F_s$ has a preferred orientation determined by the following rule: the orientation of $M$ at any point of $F_s$ coincides with the product of the orientation induced by the standard one of $S^1$ on any smooth local section of $\rho_f$ with the preferred orientation of $F_s$, in that order. This makes $\rho_f$ into an oriented locally trivial fibre bundle. In what follows, we will consider $F = F_*$ endowed with this preferred orientation, for the fixed base point $* \in S^1$.

Two open books $f: M \to B^2$ and $f': M' \to B^2$ are said to be *fibred equivalent* if there are orientation preserving diffeomorphisms $\varphi: B^2 \to B^2$ and $\widetilde{\varphi}: M \to M'$ such that $\varphi \circ f = f' \circ \widetilde{\varphi}$.

By the *monodromy* of an open book $f: M \to B^2$ with binding $L \subset M$ and page $F$, we mean the mapping class $\gamma_f = \omega_f(\alpha) \in \mathcal{M}_+(F)$, with $\omega_f: \pi_1(S^1) \to \mathcal{M}_+(F)$ the monodromy homomorphism of the $F$-bundle $\rho_{f|M-\mathrm{Int}\,T}\colon M - \mathrm{Int}\,T \to S^1$ and $\alpha \in \pi_1(S^1)$ the usual counterclockwise generator. The monodromy $\gamma_f$ uniquely determines $\rho_f$ and hence $f$ up to fibred equivalence.

Once an identification $F \cong F_{g,b}$ is chosen, we can think of $\gamma_f$ as an element of $\mathcal{M}_{g,b} = \mathcal{M}_+(F_{g,b}) \cong \mathcal{M}_+(F)$. Of course, this is only defined up to conjugation in $\mathcal{M}_{g,b}$, depending of the identification $F \cong F_{g,b}$. Actually, two open books $f\colon M \to B^2$ and $f'\colon M' \to B^2$ are fibred equivalent if and only if they have diffeomorphic pages $F \cong F' \cong F_{g,b}$ and conjugate monodromies $\gamma_f$ and $\gamma_{f'}$ in $\mathcal{M}_{g,b}$.

Given any $\gamma \in \mathcal{M}_{g,b}$, we can construct an open book $f_\gamma\colon M_\gamma \to B^2$ with page $F \cong F_{g,b}$ and monodromy $\gamma$ as follows. Let $T(\gamma) = F_{g,b} \times [0,1]/((\gamma(x),0) \sim (x,1)\ \forall x \in F_{g,b})$ be the mapping torus of (any representative of) $\gamma$. Since $\gamma$ restricts to the identity of $\mathrm{Bd}\,F_{g,b}$, there is a canonical identification $\eta\colon \mathrm{Bd}\,F_{g,b} \times S^1 \to T(\gamma_{|\,\mathrm{Bd}\,F_{g,b}}) \subset T(\gamma)$. Then, we put $M_\gamma = T(\gamma) \cup_\eta \mathrm{Bd}\,F_{g,b} \times B^2$ and define $f_\gamma\colon M_\gamma \to B^2$ to coincide with the canonical projection $\pi\colon T(\gamma) \to S^1 \cong [0,1]/(0 \sim 1)$ on $T(\gamma)$ and with the projection on the second factor on $\mathrm{Bd}\,F_{g,b} \times B^2$.

It is clear from the definitions that any Lefschetz fibration $f\colon W \to B^2$ restricts to an open book $\partial f = f_|\colon \mathrm{Bd}\,W \to B^2$ on the boundary of $W$. The page of $\partial f$ is the regular fibre $F$ of $f$ and its monodromy homomorphism is $\omega_{\partial f} = \omega_f \circ i_*\colon \pi_1(S^1) \to \mathcal{M}_{g,b}$, where $i_*$ is the homomorphism induced by the inclusion $i\colon S^1 \subset B^2 - A$, with $A$ the set of singular values of $f$. Hence, if $(\gamma_1, \gamma_2, \ldots, \gamma_n)$ is any mapping monodromy sequence for $f$, the monodromy of $\partial f$ is given by the product $\gamma_{\partial f} = \gamma_1 \gamma_2 \cdots \gamma_n$ (usually called the *total monodromy* of $f$).

Conversely, any open book $f\colon M \to B^2$ can be easily seen to be fibred equivalent to the boundary restriction $\partial \widetilde{f}\colon \mathrm{Bd}\,W \to B^2$ of an allowable Lefschetz fibration $\widetilde{f}\colon W \to B^2$, which we call a *filling* of $f$. In fact, $\mathcal{M}_{g,b}$ is known to be generated by the Dehn twists along homologically non-trivial cycles in $F_{g,b}$, and any factorization $\gamma_f = \gamma_1 \gamma_2 \cdots \gamma_n \in \mathcal{M}_{g,b}$ of the monodromy of $f$ into such Dehn twists gives rise to a mapping monodromy sequence $(\gamma_1, \gamma_2, \ldots, \gamma_n)$ representing a filling of $f$. Of course, different factorizations give rise to possibly inequivalent different fillings.

In particular, we consider the standard open book $\partial \pi\colon S^3 \cong \mathrm{Bd}(B^2 \times B^2) \to B^2$ on $S^3$ as the boundary of the trivial product fibration $\pi\colon B^4 \cong B^2 \times B^2 \to B^2$ given by the projection onto the first factor.

If $f\colon W \to B^2$ is an allowable Lefschetz fibration and $f = \pi \circ p$ is any factorization with $p\colon W \to B^2 \times B^2$ a simple covering branched over a braided surface $S \subset B^2 \times B^2$ as in Proposition 11, then the boundary open book $\partial f$ admits an analogous factorization $\partial f = \partial \pi \circ p_|$, where $p_|\colon \mathrm{Bd}\,W \to S^3 \cong \mathrm{Bd}(B^2 \times B^2)$ is a simple covering branched over the closed braid $\mathrm{Bd}\,S \subset S^1 \times B^2$. By the existence of fillings, any open book admits such a factorization, hence it can be represented as the lifting of the standard open book $\partial \pi$ with respect to a simple covering of $S^3$ branched over a closed braid.

It is well known since Alexander [1] that any closed oriented three-manifold admits an open book decomposition (Gonzáles-Acuña [10] and Myers [18] independently proved that this can always be assumed to have connected binding).

Harer in [12] proved that any two open book decompositions of the same three-manifold are related, up to fibred equivalence, by a sequence of Hopf band (de)plumbing and double twistings. The Hopf band (de)plumbing is essentially the restriction to the boundary of the Hopf (de)stabilization of Lefschetz fibrations described in Section 7. While the double twisting is a more involved modification defined in terms of surgery as follows.

Given an open book decomposition of $M$, consider two pages $F_1$ and $F_2$ of it and two cycles $c_1 \subset F_1$ and $c_2 \subset F_2$ which bound an embedded annulus $A \subset M$. Let $c'_i \subset F_i$ and $c''_i \subset A$ the framings induced on $c_i$ by $F_i$ and $A$, respectively, and assume that $c''_i + (-1)^i = c'_i + \varepsilon_i$, for some arbitrary independent choices of $\varepsilon_i = \pm 1$. Then, surgering $M$ along $c_1$ and $c_2$ with the opposite framings $c''_1 - 1$ and $c''_2 + 1$ does not change the manifold $M$, while it changes the original monodromy $\gamma \in \mathcal{M}_{g,b}$ of the open book into the composition $\gamma \gamma_1^{\varepsilon_1} \gamma_2^{\varepsilon_2} \in \mathcal{M}_{g,b}$, where $\gamma_i$ is the positive Dehn twists of $F_{g,b}$ along $c_i \subset F_i \cong F_{g,b}$ for a suitable identification $F_i \cong F_{g,b}$.

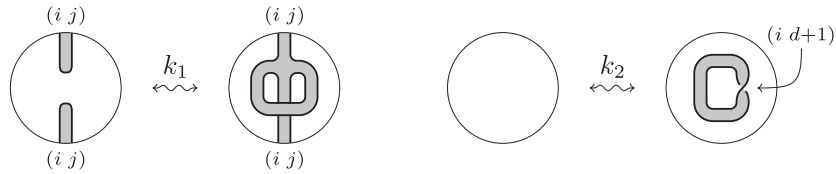

Figure 63. *Kirby calculus moves for planar diagrams.*

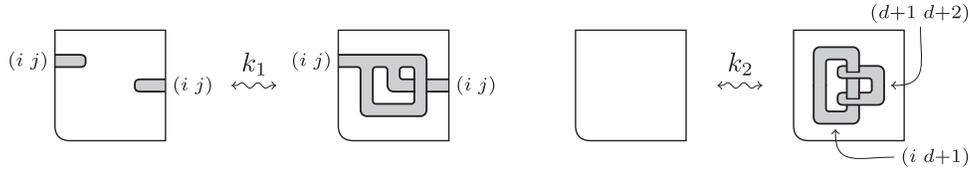

Figure 64. *Kirby calculus moves for rectangular diagrams.*

Unfortunately, the effect of a double twisting on the open book structure can be quite destructive, due to the fact that the annulus $A$ can intersect the binding and the pages in a rather arbitrary way. However, as conjectured by Harer himself in [**12**] and later proved by Giroux and Goodman in [**8**], this second move is not needed when $M$ is $S^3$ (or more generally an integral homology sphere).

Here, we will provide a different set of moves alternative to the one given by Harer, based on the results of the previous section, by looking at open books as boundaries of Lefschetz fibrations.

As the first step, we establish how to relate two Lefschetz fibrations on four-dimensional 2-handlebodies having diffeomorphic oriented boundaries. In order to do that, let us introduce the moves in Figure 63, which interpret the Kirby calculus moves in terms of labelled planar diagrams, according to next proposition. Here, on the left side of move $k_2$ we assume to have a $\Sigma_d$-labelled diagram.

PROPOSITION 17. *Two labelled ribbon surfaces in $B^4$ represent connected four-dimensional 2-handlebodies with diffeomorphic boundaries if and only if they are related by labelled 1-isotopy, (de)stabilization and the moves $c_1, c_2, k_1$ and $k_2$ in Figures 24 and 63.*

*Proof.* This is essentially [**3**, Theorem 2]. In fact: move $k_1$ coincides with move $T$ of [**3**]; move $k_2$ is equivalent to move $P_+$ of [**3**] up to stabilization and the covering move $c_3$ in Figure 25; $P_-$ of [**3**] is equivalent to the inverse of $P_+$ modulo $k_1$. □

Figure 64 shows rectangular versions of the moves above, in a suitable restricted form for the application of the braiding procedure.

LEMMA 18. *Two labelled rectangular diagrams represent connected four-dimensional 2-handlebodies with diffeomorphic boundaries if and only if they are related by rectangular (de)stabilization, the moves $r_1$ to $r_{25}$ in Figures 39–42, and the moves $k_1$ and $k_2$ in Figure 64.*

*Proof.* This is an immediate consequence of Lemma 15 and Proposition 17. □

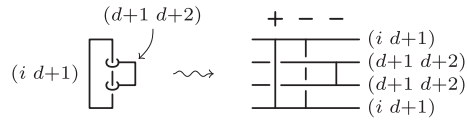

FIGURE 65. *Braiding the move $k_2$.*

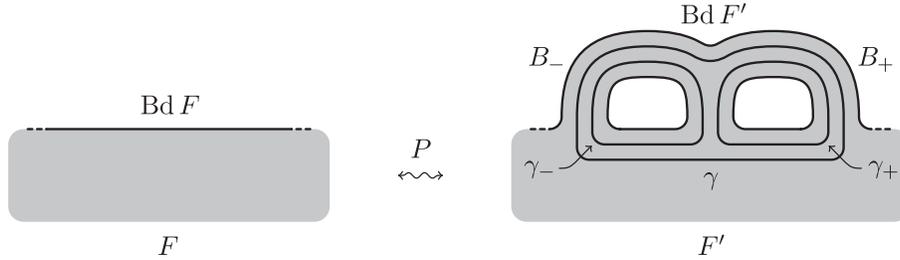

FIGURE 66. *The $P$ move.*

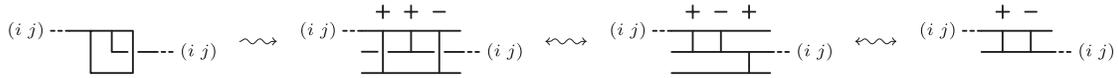

FIGURE 67. *Braiding the move $k_1$.*

Before stating our equivalence theorem for the boundaries of allowable Lefschetz fibrations, we still need to see how moves $k_1$ and $k_2$ look like once the braiding procedure is applied to them.

*P move.* We start with move $k_2$ in Figure 64. Since the right side of the move is separated from the rest of the labelled rectangular diagram, we can think of that move as adding (or deleting) the corresponding labelled braided surface shown in Figure 65 on the top of the other sheets (labelled in $\Sigma_d$).

In terms of the mapping monodromy sequence of the Lefschetz fibration $f$, this means adding (or deleting) two bands $B_-$ and $B_+$ to the regular fibre $F$, and three Dehn twists $\gamma, \gamma_-$ and $\gamma_+$ to the sequence, the first two twists negative and the third one positive, as illustrated in Figure 66. We leave the easy verification of that to the reader. Of course, being those cycles disjoint from one another and from all the other ones, it does not matter where they are located in the sequence.

We call that modification a $P$ move. Actually, another version of the $P$ move, equivalent to the inverse of it modulo the $Q$ move below, could be given with $\gamma$ a positive twist. This would correspond to inverting the half-twist in the band on the right side of the original move $k_2$ in Figure 63.

Up to $U$ move (used to get allowability), both versions of the $P$ move, with $\gamma$ negative or positive, can be thought as adding/deleting a trivial Dehn twist $\gamma$ to the monodromy sequence. This means performing a blow-up/down in $\operatorname{Int} W$, by adding/deleting a connected sum term $+\mathbb{C}P^2$ or $-\mathbb{C}P^2$, respectively, leaving the boundary open book $\partial f$ unchanged.

*Q move.* Figure 67 shows the labelled braided surface translation of the diagram of the right side of move $k_1$ in Figure 64. This differs from the analogous translation of the left side, by an extra pair of contiguous opposite Dehn twists.

We call $Q$ move the insertion (or deletion) in the mapping monodromy sequence of a Lefschetz fibration $f$ of any pair of contiguous opposite Dehn twists along a homologically non-trivial cycle. Obviously, this does not affect the total monodromy, hence it leaves the boundary open book $\partial f$ unchanged.

At this point, we are ready to conclude with the boundary equivalence theorem and its interpretation in terms of open books.

THEOREM B. *Two allowable Lefschetz fibrations over $B^2$ represent four-dimensional 2-handlebodies with diffeomorphic oriented boundaries if and only if they are related by fibred equivalence, the moves $S$ and $T$ of Section 7, and the moves $P$ and $Q$.*

*Proof.* In light of the above discussion on $P$ and $Q$ moves, this follows from Theorem A and Lemma 18. □

Now, we denote by $\partial S$, $\partial T$ and $\partial P$ the moves on open books given by the restriction to boundary of the corresponding moves $S$, $T$ and $P$ on Lefschetz fibrations. In particular, $\partial S$ coincides with the Hopf plumbing considered by Harer in [**12**], while $\partial T$ and $\partial P$ are briefly discussed below.

THEOREM C. *Two open books are supported by diffeomorphic oriented three-manifolds if and only if they are related by fibred equivalence and the moves $\partial S$, $\partial T$ and $\partial P$.*

*Proof.* This immediately follows from Theorem B, taking into account that the restriction of the $Q$ move to the boundary does not affect at all the open book structure. □

## 10. Final remarks

Compared with Harer's double twisting, our moves $\partial T$ and $\partial P$ seem preferable, since they can be completely described in terms of the open book monodromy.

Actually, the $\partial P$ move can be easily seen to be a special case of the Harer's double twisting, up to Hopf plumbing. In fact, referring to Figure 66, once the page $F$ has been stabilized to $F'$ with the new Dehn twists $\gamma_-$ and $\gamma_+$ by two $\partial S$ moves, the cycle $\gamma$ spans a disc $D$ in the boundary three-manifold $M$. Moreover, $D$ and $F'$ support the same framing on $\gamma$, in such a way that the Dehn twist along $\gamma$ can be inserted in the monodromy of the open book by a double twisting with $c_1$ trivial and $c_2 = \gamma$ (cf. definition at p. 385).

A similar explicit realization of the $\partial T$ move in terms of Hopf plumbing and double twisting should likely exist, but we were not able to find it. Of course, this would lead to an alternative proof of the Harer's equivalence theorem in [**12**].

As we mentioned when defining the $T$ move, the Dehn twists involved in it could have arbitrary signs. Denote by $T_+$ the move in the case of all positive twists. Then $S_+$ and $T_+$ are moves for positive (allowable) Lefschetz fibrations with bounded fibre. By [**15**], these fibrations represent compact Stein domains with strictly pseudoconvex boundary up to orientation-preserving diffeomorphisms.

PROBLEM 19. *Do $S_+$, $T_+$ and fibred equivalence suffice to relate any two positive Lefschetz fibrations on the same four-dimensional 2-handlebody up to 2-equivalence?*

A similar question can be posed for positive open books (namely those whose monodromy admits a factorization into positive Dehn twists) by considering the boundary restrictions $\partial S_+$ and $\partial T_+$.

PROBLEM 20. Do $\partial S_+$, $\partial T_+$ and fibred equivalence suffice to relate any two positive open books on the same three-manifold?

In [**7**], Giroux proved that two open books represent the same contact three-manifold if and only if they are related by positive stabilizations and fibred equivalence (see also [**6**] for a proof). It is natural to ask if our move $\partial T_+$ preserves properties of contact structures like Stein fillability, symplectic fillability or tightness.

Finally, we note that our Theorems A, B and C are formulated up to fibred equivalence. However, when the ambient four-dimensional 2-handlebody or three-manifold is given, isotopic versions of them would be desirable, with fibred isotopy in place of fibred equivalence (in the spirit of [**8**]). In order to get such isotopic versions, we just need a careful analysis of how our moves can be realized as embedded ones. This will be the object of a forthcoming paper.

*Nikos Apostolakis*
*Department of Mathematics and Computer Science*
*Bronx Community College – CUNY*
*2155 University Avenue*
  *Bronx, NY 10453*
*USA*

nikos.ap@gmail.com

*Daniele Zuddas*
*Dipartimento di Matematica e Informatica*
*Università di Cagliari*
*Via Ospedale 72*
*Andrea Loi*
*I-09124 Cagliari*
*Italy*

d.zuddas@gmail.com

*Riccardo Piergallini*
*Scuola di Scienze e Tecnologie*
*Università di Camerino*
*via Madonna delle Carceri*
*I-62032 Camerino*
*Italy*

riccardo.piergallini@unicam.it